\documentclass[12pt,a4paper]{article}
\usepackage{pdflscape}
\usepackage{amsmath}
\usepackage{amssymb}
\usepackage{latexsym}
\usepackage{srcltx}

\usepackage{subcaption}
\usepackage{tikz}
\usepackage{anyfontsize}
\usepackage{tabularx}

\usepackage{multirow}
\usepackage{array}

\usepackage{amssymb, amsmath, amsfonts, mathtools}
\usepackage{esint}

\usepackage{breqn}
\usepackage{float}
\usepackage{caption}
\usepackage{subcaption}
\usepackage{hhline}
\usepackage{multirow}
\usepackage{rotating}
\usepackage{fixltx2e}

\newcommand{\co}{{\mathbb C}}

\newcommand{\re}{{\mathbb R}}
\newcommand{\n}{{\mathbb N}}

\newcommand{\z}{{\mathbb Z}}
\newcommand{\cA}{{\mathcal{A}}}

\newcommand{\cG}{{\mathcal{G}}}
\newcommand{\cI}{{\mathcal{I}}}

\newcommand{\cM}{{\mathcal{M}}}

\newcommand{\cB}{{\mathcal{B}}}
\newcommand{\cT}{{\mathcal{T}}}
\newcommand{\cH}{{\mathcal{H}}}

\newcommand{\cR}{{\mathcal{R}}}

\newcommand{\bc}{{\boldsymbol c}}
\newcommand{\bd}{{\boldsymbol d}}

\newcommand{\bal}{{\boldsymbol \alpha}}
\newcommand{\bs}{{\boldsymbol s}}

\newcommand{\bT}{{\boldsymbol{T}}}

\newtheorem{theorem}{Theorem}
\newtheorem{prop}{Proposition}
\newtheorem{lemma}{Lemma}
\newtheorem{cor}{Corollary}
\newtheorem{remark}{Remark}
\newtheorem{ex}{Example}
\newtheorem{defi}{Definition}
\newtheorem{conj}{Conjecture}

\date{}

\author{
Vladimir Yu.~Protasov
\thanks{University of L'Aquila (Italy), 
Moscow State University (Russia),  {e-mail: \tt\small
v-protassov@yandex.ru}}}

\title{Surface dimension,  tiles, and synchronising automata
\thanks{
The research is supported by the Russian
Foundation for Basic Research, projects  no. 19-04-01227 and 20-01-00469.
}}

\begin{document}
\maketitle

\begin{abstract}

We study the surface regularity of compact sets~$G \subset \re^n$ 
which is equal to the supremum of numbers~$s\ge 0$ such that 
the measure of the set $G_{\varepsilon}\setminus G$
does not exceed~$\, C\, \varepsilon^{\,s}\, , \, \varepsilon > 0$, where 
$G_{\varepsilon}$ denotes the $\varepsilon$-neighbourhood of~$G$. The surface dimension 
is by definition the difference between~$n$ and the surface regularity. 
Those values provide a natural characterisation of regularity for sets 
of positive measure. We show that for self-affine attractors and tiles 
those characteristics are explicitly computable and find them for some popular tiles. 
This, in particular, gives a refined regularity  scale for the multivariate 
Haar wavelets. The classification  of attractors of the highest possible regularity is addressed. 
The relation between the surface regularity and the H\"older regularity 
of multivariate refinable functions and wavelets is found. Finally, the surface regularity 
is applied to the theory of synchronising automata, where it corresponds to the 
 concept of parameter of synchronisation.

\bigskip

\noindent \textbf{Keywords:} {\em Regularity, dimension, surface area, Minkowski content, self-affine attractors, tiles,  spectral radius, multivariate Haar wavelets, finite deterministic automata, reset word, synchronisation}
\smallskip

\begin{flushright}
\noindent  \textbf{AMS 2010} {\em subject
classification:  28A75,  39A99, 11K55,  68Q45}

\end{flushright}

\end{abstract}
\bigskip

\begin{center}
\textbf{1. Introduction}
\end{center}
\bigskip 

The  well-known Minkowski -– Steiner formula defines the area, i.e.,  the $(n-1)$-dimensional 
volume, of the surface of a compact set~$G \subset \re^n$ as the 
lower limit for the ratio $\frac{|G_{\varepsilon}| - |G|}{\varepsilon}$ as 
$\varepsilon \to +0$, where $G_{\varepsilon}$ is  the $\varepsilon$-neighbourhood 
of $G$ and $|X|$ denotes the Lebesgue measure of the set~$X$ (see, for instance, 
\cite{F}).  
For sets with sufficiently regular surfaces, this limit is finite. This is the case, for example, if $G$ is convex. If this limit is infinite, then a natural characterisation of regularity of the surface is the supremum of 
$s\ge 0$ such that $|G_{\varepsilon}| - |G|$ does not exceed  $C \, \varepsilon^{\,s}$ for all~$\varepsilon > 0$. This  is, in a sense, analogous 
to the H\"older exponent  of a function while the surface area plays the role 
of Lipschitz constant. In this paper we show that for self-affine tiles 
and attractors, this characteristic is computable and gives a natural 
scale of regularity for those sets. It is related to the H\"older regularity of Haar wavelets in~$\re^n$. Moreover, this characteristic can be applied 
in the study of 
synchronising automata, where it corresponds to their ``rate of synchronisation''.

We use the following notation: $\cB(x, r)$ is 
the Euclidean ball of radius $r>0$ centered at a point $x \in \re^n$; 
$A+B = \{a+b\, | \, a\in A, \, b\in B\}$ is the Minkowski sum of sets 
$A$ and $B$; $\, G_{\,\varepsilon}\, = \, G \, + \, \cB(0, \varepsilon)$ is the 
$\varepsilon$-neighbourhood  of a set~$G$. 
 \begin{defi}\label{d.s} 
The {\em surface regularity} of a compact set~$G \subset \re^n$ is 
$$
\bs(G) \ = \ \sup\, \bigl\{ s \, \ge \, 0 \ \bigr| \quad  
|G_{\varepsilon}| \, - \, |G| \ \le \ C\, \varepsilon^{\,s} \quad 
\forall \varepsilon > 0 \bigr\}. 
$$
The {\em surface dimension} of~$G$ is $\, \bd \, = \, n\, - \, \bs(G)$. 
\end{defi}
The surface dimension~$\bd$ of compact sets in~$\re^n$ 
can take all values from $0$ to $n$. 
The case of integer  $\bd$ corresponds to the upper Minkowski 
content, see~\cite{AHK, ACV, KP, R}. However, for the sets~$G$ of positive measure, we always 
have $\bs(G) \le 1$ (Corollary~\ref{c.05} in the next section), 
and therefore, $\bd(G) \in [n-1, n]$.  The case~$\bs = 1$ (i.e., $\bd = n-1$) characterises 
 sets with ``regular'' surfaces. For example, if $G$ is a union of finitely many convex sets or sets with piecewise-smooth boundary, then~$\bs=1$. 
One can say that $\bs(G)$ measures the regularity of the boundary of~$G$.
On the other hand, it has no relation to the dimension  of the topological boundary. This can be shown  by simple examples. Consider the following ``quasi-Cantor'' 
set~$G\subset \re$: 
    take a unit segment, remove the open interval of length $2^{-2^1}$
    from the middle; in each of the two remaining segments remove
    the interval of length~$2^{-2^2}$ from the middle, etc. 
    In $k$th iteration we have $2^{k}$ equal segments and from each of them we remove an interval of length ~$2^{-2^k}$ from the middle. The limit compact set~$G$
has a positive measure. It is easily shown  that~$\bs = 1$ and 
hence $\bd = 0$. On the other hand, the boundary of~$G$ coincides with~$G$
and hence the Hausdorff dimension of its  boundary is one and so it is not equal to~$\bd$.

  An advantage of the surface regilarity and of the surface dimension is that 
  they are both metric invariants of compact sets, i.e., these charactersistics are 
  invariant with respect to bi-Lipschitz maps
  (Lipschitz maps with Lipschitz inverse). Hence, they provide characteristics of compact sets invariant under $C^1$-diffeomorphisms. In contrast to the 
 topological or  Hausdorff dimension it can distinct sets of positive measure, whose 
 dimension is the same as the dimension of the entire space. That is why our main 
 interest is in the sets of positive Lebesgue measure. 
 In what follows we assume that $|G| > 0$.

 For sets of positive measure, a 
 characteristic similar to the surface regularity is provided by the $L_1$-H\"older regularity of the indicator function. However,  
 this charactersistic, in contrast to the surface regularity, is not 
 bi-Lipschitz invariant. 
 We show that the surface regularity  does not exceed 
 the H\"older regularity and can be strictly smaller even for tiles (Theorem~\ref{th.10}).
 Then in Theorem~\ref{th.20} we establish a condition for a compact set 
 that ensures that its surface and H\"older regularities coincide. 
 In Section~3 we apply that result to self-affine attractors and tiles, which play an important role in construction of Haar and other wavelet systems in~$\re^n$.
 In Section~5 we obtain formulas for the $L_p$-H\"older regularity of attractors and tiles. This, in particular, makes it possible to compute the 
 $L_p$-exponents of multivariate Haar wavelets and to range them 
 by their regularity. 
 In case of isotropic dilation matrix, those formulas can compute the surface regularities and surface dimensions. The computation of all those  
 characteristics are reduced to finding the Perron eigenvalue of a special matrix. 
 In Section~6 we compute surface dimensions of some popular 
 self-affine tiles. Then we address the problem of characterising 
 the self-affine attractors and tiles with the highest surface regularity. 
 We make a conjecture that the only self-affine attractor with the 
 surface regularity~$\bs = 1$ is a parallelepiped. So, the parallelepiped is the only regular attractor. In Section~7 this conjecture is 
 proved for dimension~$n=1$. Finally, we apply the surface regularity in the study of 
 finite deterministic automata and establish a relation between 
 the surface regularity and the rate of synchronisation (Section~8). 
 
 The following notation will be used: $|X|$ is the Lebesgue measure of a set~$X$ or the 
 cardinality of a finite set~$X$, depending on the context; $L_p$ is the 
 standard functional space with the norm $\|f\|_p \, = \, 
 \bigl(\int |f|^p dt\bigr)^{1/p}$. 
 We use the standard notation~$A^* = \bar A^T$ for the adjoint matrix to a matrix~$A$;   the spectral radius of~$A$, i.e., the biggest modulus of eigenvalues, 
 is denoted my~$\rho(A)$.

\bigskip

\begin{center}
\textbf{2. The surface dimension and the H\"older regularity}
\end{center}
\bigskip 

The definition of the surface regularity is similar to 
the H\"older regularity of the characteristic function $\chi_G(x)$ in the 
space~$L_1(\re^n)$. As usual, $\chi_G(x) \, = \, 
1$ if $x \in G$ and $\chi_G(x) \, = \, 0$ otherwise. 
The $L_p$ H\"older regularity of a function $f\in L_p(\re^n)$
is defined as 
$$
\bal_p(f) \ = \ \sup\, \bigl\{ \alpha \ge 0 \ \bigr| \ 
\|f(\cdot + h) - f(\cdot)\|_{p}\, \le \, C\, \|h\|^{\alpha} \ 
\forall h \in \re^n\, \bigr\}
$$
For a characteristic function of a compact set~$G$, we denote shortly 
$\bal_p(\chi_G) = \bal_p(G)$ and call this value 
the H\"older $L_p$-regularity of~$G$. 
The measure of the difference $G_{\varepsilon}\setminus G$ is the 
$L_1$-norm of the function $\chi_{G_{\varepsilon}\setminus G}$. 
Hence it is quite expected that $\bs(G)$ can be related 
to $\bal_1(G)$. In what follows we omit the index $1$ meaning that 
always $p=1$ if the the converse is not stated. Moreover, often 
we omit the set $G$ from the notation. Thus $\bal_1(G) = \bal$. 

The following proposition shows that 
for every compact set, the H\"older regularity majorates the surface regularity.      

\begin{prop}\label{p.10}
For every compact set in $\re^n$ of positive measure, we have $\bs \le \bal$. 
 \end{prop}  
\noindent {\tt Proof.} Let $h\in \re^n$ be an arbitrary vector 
of length~$\|h\| < \varepsilon$. 
Since $G+h \, \subset \, G+\cB(0,\varepsilon)$, we see that 
the measure of the set $(G+h)\setminus G$ does not 
exceed the measure of $G_{\varepsilon}\setminus G$. 
Similarly, the measure of $(G-h)\setminus G$ does not 
exceed the same measure of $G_{\varepsilon}\setminus G$.
Therefore, $\bigl\|\chi_{G}(\cdot + h) - \chi_{G} \bigr\|_1 \, \le \, 
2\, |G_{\varepsilon}\setminus G|$. Computing logarithms of 
both parts and dividing by 
$\log \, \frac{1}{\varepsilon}$, we conclude the proof.

{\hfill $\Box$}
\smallskip

Since the H\"older exponent never exceeds~one, we obtain 
\begin{cor}\label{c.05}
For a compact set of positive measure, $\bs \le 1$
and respectively $n-1 \le \bd \le n$. 
\end{cor}
In the sequel we always consider sets of positive measure, i.e., 
assume that $|G| > 0$. 
There are examples when $\bs \ne \bal$. 
Moreover, even for tiles on~$\re$, it can happen that~$\bs < \bal$. 
A compact set in~$\re$ is called a tile if its integer shifts
cover~$\re$ with intersections of zero measure.  

 \begin{theorem}\label{th.10}
There is a tile in~$\re$ for which $\bal = 1$ and 
 $\bs = \frac{1}{3}$. 
 \end{theorem}  
\noindent {\tt Proof.}  
First we construct a compact set~$G\subset \re$ with this property and then 
make a tile from it. 
Consider a sequence~$x_1, x_2, \ldots $, where 
 $x_k = \sum_{m=1}^k \frac{1}{m^2}$. Define  the set $G$ as a union of segments 
 $\bigl[x_k\, , \, x_k+2^{-k-2}\bigr]\, , \, k \in \n$,  plus the limit point 
 $x_{\infty} = \sum_{m=1}^{\infty} \frac{1}{m^2} \, = \, \frac{\pi^2}{6}$.

Let us first compute $\bs (G)$. 
Take arbitrary small~$\varepsilon > 0$ and denote by $N$
the minimal natural number such that for all $k \ge N$
the distance between points $x_k$ and $x_{k+1}$ is less than $2\varepsilon$. 
Thus, $N$ is the smallest natural solution of inequality 
$\frac{1}{k^2} - \frac{1}{(k+1)^2} + 2^{-k-2}\, < \, 2\varepsilon$. 
It is shown easily that  $N \, \sim \, \varepsilon^{-1/3}$ as $\varepsilon \to 0$. 
The enlarged set $G_{\varepsilon} = G + [-\varepsilon, \varepsilon]$
contains $N-1$  segments:  
$$
\left[x_k-\varepsilon\, , \, x_k+2^{-k-2}+\varepsilon\right], \quad k = 1, \ldots , N-1\, ;
$$ 
and one big segment $[x_N - \varepsilon, x_{\infty} + \varepsilon]$
formed by all other segments for $k \ge N$. 
The total length of those $N$ segments 
is $\, 2\varepsilon \, + \, \sum_{k=N}^{\infty} \frac{1}{k^2}$
(the big segment) plus $\, \sum_{k=1}^{N-1} (2^{-k-2}+2\varepsilon)$
(the remaining $N-1$ segments). Thus, 
$$
|G_{\varepsilon}| \quad = \quad 
2\varepsilon \, + \, \sum_{k=N}^{\infty} \frac{1}{k^2} \, + \, 
 \sum_{k=1}^{N-1} (2^{-k-2}+2\varepsilon) \quad = \quad 
2\, N\, \varepsilon \ + \ 
 \frac14 \, - \, 2^{-N-2}\ + \ 
\sum_{k=N}^{\infty} \frac{1}{k^2}\, .
$$ 
On the other hand, $|G| = \frac{1}{4}$. Hence 
$$
|G_{\varepsilon}| \ - \ |G| \quad =\quad  
2\, N\, \varepsilon \ - \ 2^{-N-2} \ +\ \sum_{k=N}^{\infty} \frac{1}{k^2}\, .
$$
Since $N \asymp \varepsilon^{\, -1/3}$ and 
$\sum_{k=N}^{\infty} \frac{1}{k^2} \, \asymp \, \frac{1}{N}$, we see that 
the value $|G_{\varepsilon}|- |G|$
is asymptotically 
equivalent to $\, 2\, \varepsilon^{\, 2/3}\, + \, C \, \varepsilon^{1/3}$, 
where $C$ is a constant. We see that $\, |G_{\varepsilon}|- |G| \, \asymp \, 
 C \, \varepsilon^{1/3}$ as $\, \varepsilon \to 0$, and 
therefore, $\bs = \frac{1}{3}$. 
 
Now let us  show that $\bal(G) = 1$. If $k \le  \log_2 \frac{1}{\varepsilon}$, 
i.e.,  $\, 2^{-k} > \varepsilon$, then the
 $k$th segment $[x_k, x_{k}+ 2^{-k}]$ intersects 
its copy shifted by $\varepsilon$. Therefore, the 
length of the symmetric difference of this segment with its copy is 
equal to $2\varepsilon$. The number of those segments does not exceed 
$\log_2 \frac{1}{\varepsilon}$. Hence, the total length of those 
symmetric differences is at most $\, 2\varepsilon\, 
\log_2 \frac{1}{\varepsilon}$. 

The total length of the remaining segments of the set~$G$
is $\sum_{k >  \log_2 \frac{1}{\varepsilon}}2^{-k}$, which is less then 
$2^{1-\log_2 \frac{1}{\varepsilon}} = 2 \varepsilon$. 
Therefore, the symmetric difference of this set with its 
shift to $\varepsilon$ has the length less that $4 \varepsilon$. 

Summing over these two sets we have $\|\chi_{G}(\cdot) - \chi_G(\cdot + \varepsilon)\|_1 \, <\,  2\varepsilon\, \bigl( 2 \, + \, 
\log_2 \frac{1}{\varepsilon}\bigr)$. Therefore, $\bal \ge 1$. 
Since $\bal$ cannot be bigger than one, we conclude that $\bal = 1$. 

Thus, a compact set~$G$ with $\bal = 1$ and 
 $\bs = \frac{1}{3}$ is constructed. But this is not a tile. 
 To make a tile form~$G$ we take the unit segment $[0,1]$, unify it 
with the set $G$ and  subtract the set $G-1$ from it.  
 Since $G \subset [0,2)$, we see that 
the obtained set $[0,1]\cup G \setminus (G-1)$ is a tile 
with the same parameters $\bs$ and $\bal$.

{\hfill $\Box$}
\smallskip

\begin{remark}\label{r.5}{\em 
The fact that $\bs \le \bal$ means that the characteristic 
$\bs(G)$ provides a more refined analysis of a set~$G$ than $\bal(G)$
and can distinct sets with the identical exponent~$\bal$. For example, 
the set~$G$ constructed in the proof of Theorem~\ref{th.10} 
has the maximal H\"older regularity~$\bal = 1$ as the segment~$[0,1]$, while 
its surface regularity is lower: $\bs  = \frac13$ for $G$
instead of $\bs = 1$ for the segment $[0,1]$. So, the H\"older regularity cannot 
distinct the set~$G$ from a segment, but the surface regularity can.  
 }
\end{remark}

A question arises what conditions of the set~$G$ would 
guarantee that $\bs(G) = \bal(G)$? Theorem~\ref{th.20}
below establishes sufficient conditions. To formulate them we need one more notation. 
Let a compact set $G \subset \re^n$ be fixed. 
For a  point $x\in \re^n$ and a number $r> 0$, 
we denote 
$$
\nu\, (x, r)\ = \ \frac{|\cB(x, r)\, \cap \, G|}{|\cB(x, r)|}.
$$
Thus, the number $\nu\, (x, r)$ shows which part of the volume of the ball~$\cB(x, r)$
is covered by~$G$. 
 \begin{theorem}\label{th.20}
If there are constants $c_1, c_2 > 0$ such that for every sufficiently small~$\varepsilon > 0$, 
the total measure of points $x$ of the set~$G_{\varepsilon}\setminus G$ 
for which $\nu\, (x, 2\varepsilon) \, \ge \, c_1$ is at least 
$c_2\, |G_{\varepsilon}\setminus G|$, then $\, \bs = \bal$. 
 \end{theorem}  
 If $x \in G_{\varepsilon}$, then the intersection of the ball
 $\cB(x, 2\varepsilon)$ with the set $G$ is at least nonempty.
 The assumption of Theorem~\ref{th.20} require that this intersection 
  has not very small volume. If this condition is fulfilled  for some significant part of points~$x$ of the set~$G_{\varepsilon}$, then $\bs(G) = \bal(G)$. 
  
  In the proof we use the following technical  
\begin{lemma}\label{l.10}
For an arbitrary compact set $G \subset \re^n$ and for every $r>0$, we have 
$$
|G_{\,2r}|\ - \   |G|  \quad  \le 
\quad  2^{n}\, \Bigl(\, |G_{r}| \ - \  |G| \, \Bigr)\, . 
$$
\end{lemma}
{\tt Proof.} Assume without loss of generality that $|G| = 1$. 
Since $G_r$ is a Minkowski sum of 
$G$ and of a ball of radius $r$, we can apply the Brunn-Minkowski inequality 
and conclude that the function $f(r) = \sqrt[n]{|G_r|}$
is concave. Hence $f(r) \, \ge \, \frac12\, \bigl( f(0) \, + \, f(2r)\, \bigr)$. 
Let $f(2r) = 1+a$. Since $f(0) = 1$, we see that $f(r) \, \ge \, 1 \, + \, 
\frac{a}{2}$. Therefore $\, 2^n\, \bigl(|G_r| \, - \, |G|\bigr)\, \ge \, 
2^n\, \bigl( \, \bigl(\frac{a}{2} + 1\bigr)^n \, - \, 1\bigr) \, = \, 
\bigl(a + 2\bigr)^n \, - \, 2^n$. Since $|G_{\, 2r}| = (a + 1)^n -  1$, it remains to establish the inequality 
$$
\bigl(a + 2\bigr)^n \, - \, 2^n \quad \ge \quad  \bigl(a + 1\bigr)^n\, - \, 1\, .   
$$
Opening the brackets we have $ \sum_{k=1}^n {n \choose k}\, 2^{n-k} a^k \, \ge 
\,  \sum_{k=1}^n {n \choose k}\, a^k \,$, which is obvious. 

{\hfill $\Box$}
\smallskip

\noindent {\tt Proof of Theorem~\ref{th.20}.} Denote $\chi(x) = \chi_G(x)$ and,  
for arbitrary $\varepsilon> 0$, consider the following integral: 
$$
I_{\varepsilon} \quad = \quad \int\limits_{(x, y) \in \re^n\times \re^n \atop \|x- y\|\, \le \, 2\,\varepsilon}
\bigl| \chi(x) - \chi(y) \bigr|\, dx \, dy
$$ 
This integral is computed over all pairs of points $(x, y)$ of the space $\re^n$
such that $\|x- y\|\, \le \, 2\varepsilon$. The function under the integral takes only two values:  zero and one. It is equal to zero if both $x$ and $y$ 
belong to $G$
or both do not. Otherwise it is equal to one. 
In particular, this function is zero, whenever  both $x, y$ are far from $G$. 
Therefore, it has a compact support, and  hence it is integrable 
over $\re^n\times \re^n$. 

We are going to prove that $\bs \ge \bal$ by computing the integral 
$I_{\varepsilon}$ in two ways. The first way is to integrate over the
variable~$x$:
$$
I_{\varepsilon} \quad = \quad 
\int_{h \in \cB_{\, 2\varepsilon}}\, \int_{x \in \re^n}\, 
\bigl| \chi(x) - \chi(x+h) \bigr|\, dx \, dh \quad = \quad 
\int_{h \in \cB_{\,2\varepsilon}}\ \bigl\| \chi(\cdot) - \chi(\cdot + h)\bigr\|_1\,  \, dh \ \le 
$$
$$
|\cB_{\,2\varepsilon}|\ 
\max_{h \in \cB_{\,2\varepsilon}}\, \bigl\| \chi(\cdot) - \chi(\cdot + h)\bigr\|_1\,  
$$
Thus, 
\begin{equation}\label{eq.1}
\max_{\|h\| \, \le \, 2\varepsilon}\ \bigl\| \chi(\cdot) - \chi(\cdot + h)\bigr\|_1
\quad \ge \quad |\cB_{\,2\varepsilon}|^{-1}\, I_{\varepsilon}\, . 
\end{equation}
Now we compute the integral $I_{\varepsilon}$ differently. 
Note that $\chi(x) - \chi(x+h) \ne 0$ if and only if precisely one of the 
points $x$ and $x+h$ is out of $G$. On the other hand, 
$\|h\| \, \le \, 2\,\varepsilon$, hence that point belongs to 
$G_{\,2\varepsilon}\setminus G$. Because of the symmetry we can always assume that 
$x \in G_{\,2\varepsilon}\setminus G$. Thus, 
$$
I_{\varepsilon} \quad = \quad 
\int_{G_{\,2\varepsilon}\setminus G}\, \left(\, \int_{\cB_{\,2\varepsilon}}\, 
\bigl| \chi(x) - \chi(x+h) \bigr|\, dh \right)\, dx \quad = \quad 
\int_{G_{\,2\varepsilon}\setminus G}\, |\cB_{\,2\varepsilon}|\ \nu(x, 2\varepsilon) \, dx\ \ge 
$$
$$
|\cB_{\,2\varepsilon}|\, \int_{G_{\varepsilon}\setminus G}\, \nu(x, 2\varepsilon)\, dx \quad 
\ge \quad |\cB_{\,2\varepsilon}|\, \cdot \, c_1\, \cdot \left| \Bigl\{ x \in G_{\varepsilon}\setminus G \ \Bigl| \ 
\nu(x, 2\varepsilon) \, \ge \, c_1 \, \Bigr\}\ \right| \ \ge \ 
c_1  \cdot \bigl| \cB_{\,2\varepsilon} \bigr| \cdot c_2  \cdot \bigl|\, G_{\varepsilon}\setminus G \, \bigr|\, . 
$$
Now invoking Lemma~\ref{l.10} we get 
 $\bigl|\, G_{\varepsilon}\setminus G \, \bigr| \, \ge 
\, 2^{-n}\, \bigl|\, G_{\,2\varepsilon }\setminus G \, \bigr|$, consequently  
$$
I_{\varepsilon}  \quad \ge \quad \ c_1\, c_2\, 2^{-n} \ \bigl| \cB_{\,2\varepsilon} \bigr|\, \cdot \, \bigl|\, G_{\,2\varepsilon }\setminus G \, \bigr|\, . 
$$
Combining this with~(\ref{eq.1}) we obtain 
$$
\max_{\|h\| \le 2\varepsilon}\, \bigl\| \chi(\cdot) - \chi(\cdot + h)\bigr\|_1
\quad \ge \quad c_1\, c_2\, 2^{-n}\, \bigl|\, G_{\,2\varepsilon }\setminus G \, \bigr|
$$
for every sufficiently small $\varepsilon > 0$. Taking now logarithm of both parts and a limit as $\varepsilon \to 0$ we complete the proof. 

{\hfill $\Box$}
\smallskip

Equality $\bs(G) = \bal(G)$ enables us to compute the surface dimension 
at least for some special classes of sets, because $\bal(G)$ can be expressed by the 
H\"older $L_2$-regularity as well as by the Sobolev regularity of the characteristic function.  One of such classes of sets is the class of self-affine attractors.  It plays an important role in many practical areas 
such as subdivision algorithms and wavelets. Moreover, for those sets 
the exponent of regularity~$\bal(G)$ can be efficiently computed. This is done 
in Section 5. Then in Section~8  we find a relation of this value 
to  synchronising  automata theory. 
\bigskip 

\newpage

\begin{center}
\textbf{3. The surface dimension of self-affine attractors}
\end{center}
\bigskip 

Self-affine attractors are compact sets in~$\re^n$ defined by an 
integer matrix and by a system of digits (integer points) associated  to that matrix. 

Let us have an integer $n\times n$ matrix $M$ which is supposed to be 
{\em expanding}, i.e., all its eigenvalues are strictly bigger than one by modulus. 
This matrix splits the integer lattice $\z^n$ into $m = |{\rm det}\, M|$  quotient classes defined by
the equivalence  $x \sim y \, \Leftrightarrow \, y - x \in M\, \z^n$.
Choosing one representative $d_i \in \z^n$ from each equivalence class,  we obtain
a {\em set of digits}~$D = \{d_i \ : i=0, \ldots, m-1\}$.  We always assume that~$0 \in D$ and naturally denote $d_0 = 0$.

\noindent For every integer point $d \in \z^n$,  we denote by $M_d$, the affine operator
$M_d\, x \, = \, Mx - d, \, x \in \re^n$. We use the notation $0.a_1a_2 \ldots =
\displaystyle \sum_{i=1}^{\infty} M^{-i}a_i$, $a_i \in D$. 
\begin{defi}\label{d.attractor}   
A self-affine attractor generated by an expanding matrix~$M$ and by a digit set~$D$ is the set 
\begin{equation}\label{eq.G}
G \ = \ G(M, D)\quad = \quad \left\{\, 0.a_1a_2 \ldots \ = \ 
\sum_{k=1}^{\infty} M^{-k}a_k  \quad : \quad
a_k \in D \right\}.
\end{equation}
\end{defi}
For any integer expanding matrix $M$ and for any digit set $D$, the self-affine attractor is a compact set with a nonempty 
interior~\cite{GroH,GroM}. 
Moreover, its Lebesgue measure is always a positive integer. 
It is seen easily that 
$\ G \, = \, \bigcup\limits_{d \in D} M_d^{\, -1} \, G$. 
Moreover, each set $M_d^{\, -1} \, G$ has measure $m^{-1}\, |G|$, 
hence the sum of measures of the sets $M_d^{\, -1} \, G\, , \ 
d\in D$, is exactly $m\, m^{-1}\, |G| \, = \, |G|$. 
Consequently, all those  sets  have intersections of zero measure. 
Thus $G$ is a disjunct (up to nill sets)  sum of equal sets 
that are affinely similar to~$G$. This justifies the terminology ``self-affine''. 
In what follows we  say shortly {\em attractor} and always mean  
self-affine attractors from Definition~\ref{d.attractor}. 
 
Thus, an attractor is the set of points form~$\re^n$ with zero integer part in their $M$-adic  expansion.  In this sense  attractors play role of the unit segment 
in $\re$, but for the space $\re^n$ equipped with the
$M$-adic system with digits from the set~$D$.

The affine similarity implies that the characteristic function 
$\varphi \, = \, \chi_{G}(x)$ of an attractor satisfies the following 
functional equation with a contraction of the argument: 
\begin{equation}\label{eq.ref0}
\varphi(x)\  = \  \sum_{d \in D} \, \varphi(Mx - d) \quad \mbox{ a.e.}\ \quad x \in \re^n
\end{equation}
This is a special case of a {\em refinement equation} (see Section~4). 
Therefore, the theory of refinement equations, which is  well developed in the 
literature, can be applied in the study of attractors. 
Integer shifts of an attractor cover the space~$\re^n$
with an integer number of layers (namely, with  $|G|$ layers). 
This means that~$\sum_{k \in \z^n} \varphi(x+k) \, \equiv \, |G|\ \mbox{ a.e.}$
If $|G| = 1$, then $G$ called a {\em self-affine tile}. 
\begin{defi}\label{d.tile}
A self-affine tile is an attractor of measure one. 
\end{defi}
Every integer expanding matrix and every set of digits 
generate an attractor, but this attractor is not always a tile. 
There is a criterion to determine whether the attractor 
generated by a matrix~$M$ and by a digit set~$D$ is a tile. 
It gives an answer in terms of eigenvalues of a certain integer matrix. 

In case $n=1$, for $M=2$, there are only two digit sets 
$D=\{0,1\}$ and $D=\{0,-1\}$ for which the generated attractor is a 
tile. Already for $M=3$, the situation is more interesting: 
both digit sets $D=\{0,1,2\}$ and $D = \{0,1,5\}$ generate tiles
and the second tile is not a segment. For $n=2,3$, every expanding 
matrix~$M$ has at least one digit set~$D$ generating a tile. 
However, for $n=4$,  there are examples of matrices for which this is not true.  
 Nevertheless, such a system of digits exists under 
quite general assumptions. For example, it exists whenever  
$|{\rm dim}\, M| > n\, $~\cite{LW}. 
This condition is indeed general taking into account that the matrix~$M$ is expanding. 
\begin{defi}\label{d.tiling}
A tiling~$\cG$ generated by  an integer expanding  matrix $M$ and by a set of digits $D$
is a collection of sets $\cG = \{k+ G\}_{k \in \z^s}$ such that

a)  the union of the sets in $\cG$ covers $\re^s$ and $\left| ( \ell+ G) \cap (k+ G)\right|\, =\, 0$, $\ell \not= k$;

b)  $\displaystyle \ G \, = \, \cup_{d \in D} M_d^{\, -1} \, G$.
\end{defi}
See~\cite{AG, CHM2, LW} for the general discussion and more references. 
The characteristic function of a tile possesses orthonormal integer shifts. 
This property makes tiles very useful in  the construction of the Multiresilutional analysis and wavelets systems on~$\re^n$. In particular, multivariate 
Haar systems are obtained directly from tiles~\cite{CHM2, GroM, KPS}.

A matrix is called {\em isotropic} if it has  equal by modulus eigenvalues 
 and no nontrivial Jordan blocks. An isotropic matrix is similar to 
 a multiple of an orthogonal matrix. Attractors and tiles generated  
 by isotropic dilation matrices~$M$ are very popular in applications  and 
 well studied 
 in the literature, see~\cite{CDM, H, W} and references therein. 
 The following theorem shows that at least for isotropic 
 dilation matrices, 
 the surface regularity of attractors is equal to the H\"older regularity.

 \begin{theorem}\label{th.30}
For every attractor with an isotropic dilation matrix, we have $\bs = \bal$. 
 \end{theorem}  
\noindent {\tt Proof.}  We need to show that an attractor generated by 
an isotropic dilation matrix satisfies assumptions of Theorem~\ref{th.20}. 
Take a small $\varepsilon$ and a point $x \in G_{\varepsilon}\setminus G$. 
By definition, there is a point $y \in G$ such that $\|x-y\| \, \le \, \varepsilon$. 
Denote $r = \rho(M)$. 
Let $k$ be the smallest number such that the diameter of the set 
$M^{-k}G$ is less that $\varepsilon$. Since $M$ is isotropic, 
we have $k \le \frac{\log \varepsilon}{\log r} + C_0$, where $C_0$
does not depend on~$\varepsilon$. The $k$th iteration of the partition 
$G = \cap_{i=0}^{m-1} M^{-1}(G+d_i)$ covers the set~$G$ with $m^k$ parts equal 
to $M^{-k}G$ each. Denote by~$G_k$ a part that contains the point~$y$. 
Since the diameter of $G_k$ is less than~$\varepsilon$
and $\|x-y\| \, \le \, \varepsilon$, we see that~$G_k \subset \cB(x, 2\varepsilon)$. On the other hand, $G_k \subset G$. Hence, $G_k$ is 
contained in the intersection of the ball~$\cB(x, 2\varepsilon)$ with $G$. 
Since $M$ is isotropic, we see that 
$$
|G_k| \ = \ m^{-k} |G| \ = \ r^{-nk}|G| \ \ge \ 
{\varepsilon}^{\,n}r^{-nC_0}|G|
$$
(we used the inequality $k \le \frac{\log r}{\log \varepsilon } + C_0$). 
Thus, the intersection of the ball $\cB(x, \,2\varepsilon)$ with the 
set $G$ has the measure at least $C\varepsilon^{d}$, where $C$ is a constant 
not depending on~$\varepsilon$. On the other hand, $|\cB(x, 2\varepsilon)| = 
\varepsilon^n |\cB|$. Hence, the ratio of measures of the 
intersection to the measure of the ball is at least $\frac{C}{|\cB|}$. 
This is true for all points $x\in \bigl(G_{\varepsilon}\setminus G\bigr)$. 
Therefore, the assumptions of Theorem~\ref{th.20} are satisfied with the 
constants $\, c_1 = \frac{C}{|\cB|}, \, c_2 = 1$.

{\hfill $\Box$}
\smallskip

We believe that the assumption of isotropic dilation matrix~$M$
can be omitted in Theorem~\ref{th.30} and that actually 
$\bs = \bal$ for an arbitrary  attractor. 
\begin{conj}\label{conj.5}
Theorem~\ref{th.30} holds for arbitrary dilation matrix. 
\end{conj}

In the next section we show that the H\"older regularity of attractors 
can be efficiently computed. It can be expressed with the Perron eigenvalue of a
special  matrix. In case of an isotropic matrix~$M$ this
will give the values of the surface regularity and surface dimension
of the attractors.  

\bigskip 

\begin{center} 
\textbf{4. The $L_p$-regularity of multivariate refinable functions}
\end{center}
\bigskip 

In this section we consider the $L_p$-regularity of solutions  of general 
refinement equations. We provide a method that allows, at least theoretically, 
to find the $L_p$-H\"older exponent of wavelets and of the limit functions of subdivision schemes on~$\re^n$, see~\cite{CP, KPS}. Then, in Section~5, we apply the obtained results to 
the special refinement equations~(\ref{eq.ref0}) for 
characteristic functions of attractors. As we will see, in  that case 
the H\"older regularity in~$L_1$ can be found within polynomial time.
This, in particular, gives formulas the for $L_p$-H\"older regularity 
of Haar functions. 
This also makes it possible to compute the surface dimension of tiles provided~$M$ 
is isotropic.  
\bigskip 

\begin{center}
\textbf{4.1 Refinement equation}
\end{center}
\bigskip 

 {\em Refinement equation}  is a
functional equation of the type  
\begin{equation}\label{eq.ref}
\varphi(x)\ = \ \sum_{k \in \z^n} \, c_k \, \varphi \, (Mx \, - \, k), \quad  x \in \re^n,
\end{equation}
with a compactly supported set of coefficients $c_k \in \co$ 
(i.e., $c_k = 0$ for all but finitely many~$k$) and with a general integer expanding dilation 
matrix~$M $. 
The set  of coefficients~$\bc \, = \, \{c_k, \, k \in \z^n\}$ 
is called a {\em mask} of the equation. The theory of refinement equations is 
well developed in the literature due to their crucial role in the construction of wavelets~\cite{CHM2, D, W}, in the study of  subdivision schemes for approximating functions and for curves and surfaces design~\cite{CDM, CM, J}, in some problems of combinatorics, number theory, and probability~(see~\cite{CDM, P06} and references therein). 
The characteristic function of an attractor~$\varphi = \chi_{G}$
is a solution of refinement equation with 
$c_k = 1$ whenever $k\in D$ and $c_k = 0$ otherwise~(\ref{eq.ref0}).

A compactly supported function  $\varphi \in L_1(\re^n)$ satisfying equation~(\ref{eq.ref}) is
called {\em a refinable function}. It is well known that if such a solution 
exists, then it is unique up to normalization. Moreover, if 
 $\int_{\re^n} \varphi (x)\, dx \ne 0$, then 
$\displaystyle \sum_{k \in \z^n} c_k = m$. We will focus on this case 
as in the most of literature. Under this assumption, the refinement equation always possesses 
a unique up to multiplication by a constant solution~$\varphi$ in the space of tempered distributions. This solution is compactly supported~\cite{NPS}. 
For the special case~(\ref{eq.ref0}) we have $c_j = 1$
if $j\in D$, otherwise $c_j = 0$, and the (unique!)  solution~$\varphi$ is a characteristic 
function of the attractor~$G = G(M, D)$. Furthermore, as in the most of literature
we assume that the refinement equations satisfy the {\em sum rules}: 
\begin{equation}\label{eq.sumrules}
\sum_{k\in \z^n}\, c_{\, Mk - d} \quad = \quad 1\, , \qquad d \in D\, . 
\end{equation} 
The equation for attractors~(\ref{eq.ref0}) always satisfies it  because 
the set of coefficient $\bigl\{ c_{\, Mk - d}, , \, k\in \z \bigr\}$ 
consists of zeros except for one coefficient being equal to one. 

\bigskip 

\begin{center}
\textbf{4.2 The basic tile}
\end{center}
\bigskip 

There are several methods to analyse regularity 
of solutions of  refinement equations. Some of them such as the matrix 
method can find the precise values of the H\"older exponent. The main idea  is to pass from the refinement equation on $\re^n$ to 
an equation on a vector function defined on some basic tile. So, the matrix  method requires an auxiliary tile~$Q = Q(M, \Delta)$  generated by the same matrix~$M$ and 
by some set of digits  $\Delta = \{\delta_0, \ldots , \delta_{m-1}\}$. 
An arbitrary tile generated by the  
matrix~$M$ can play the role of a basic tile. 

\bigskip 

\begin{center}
\textbf{4.3 Invariant subsets of~$\z^n$}
\end{center}
\bigskip 

The first step to realize the matrix method is to choose a special finite 
subset of~$\z^n$.   
Let us have a refinement equation~(\ref{eq.ref}) with a 
mask~$\bc \, = \, \{c_k \, , \, k \in \z^n\}$.  
Consider a map $\eta: \, 2^{\z^n}\to 2^{\z^n}$ that to every set of integers $X \subset \z^n$ 
associates the set $M^{-1}\bigl(X + {\rm supp} \, \bc \, - \, \Delta\bigr)\, \cap \, \z^n$.  Let us recall that ${\rm supp} \, \bc \, = \, 
\{k \in \z^n \, , \, c_k \ne 0\}$. Since $\bc$ is compactly supported, 
$|{\rm supp} \, \bc|\, < \, \infty$.  
\begin{defi}\label{d.inv-set}
Let a digit set $\Delta$ and a compactly supported mask $\bc \, = \, \{c_k\}_{k \in \z}$ be given.  A finite set $S \subset \z^n$ is called invariant if 
$\eta\,  S\,  \subset \, S$. 
\end{defi}
\begin{ex}\label{ex.inv-set1}
{\em Consider the real line with the dilation $M=2$ and digits $\Delta =\{0,1\}$. 
 Then for ${\rm supp}\, \bc \, = \, \{0, N\}$, the set
 $S = \{0, \ldots , N-1\}$ is invariant. Indeed, $S \, + \, {\rm supp} \, \bc \, = \, 
 \{0, \ldots , 2N-1\}$, hence $S\,  + \, {\rm supp} \, \bc \, - \, \Delta \, = \, 
 \{-1, 0, \ldots , 2N-1\}$. All integers in the set 
$\frac12 \, \{-1, 0, \ldots , 2N-1\}$  are 
$\{0, \ldots , N-1\}$, hence $\eta S \, = \,  S$. 
Every segment of integers that contains $S$ is also an invariant set.  The same holds 
for every mask with  support that contains numbers $0, N$ and 
some integers (may be all) between them. 

For the support ${\rm supp}\, \bc \, = \, \{0, 7\}$ and $\Delta=\{0,1\}$, but with the 
dilation $M = -2$, the set $S = \{0, \ldots , N-1\}$ is nor longer  
invariant, but the set $S = \{-4, -2, -1, 0, 1\}$ is. 

For the real line with the dilation $M = 3$, digits $\Delta = \{0,1,2\}$
and ${\rm supp}\, \bc \, = \, \{0, 1,5\}$, the set $S = \{0,1,2\}$
is invariant. 
}
\end{ex}
For a given mask $\bc$, we consider the following set:  
\begin{equation}\label{eq.K}
 \Gamma \ = \  \Gamma(M, {\rm supp} \, \bc) \ =\ \bigl\{\ x \in \re^n \, \bigl| \  x\ =\ \sum_{j=1}^\infty \
 M^{-j} \gamma_j, \ \  \gamma_j \in {\rm supp} \, \bc\,  \bigl\}, 
\end{equation}
If ${\rm supp} \, \bc$ is a digit set for $M$, then $\Gamma$ is an attractor. 
The set $\Gamma$ will be referred to as {\em support set} of the refinement equation.
In general, this set may not coincide with ${\rm supp} \, \varphi$. 
However, we always have ${\rm supp} \, \varphi \, \subset \, \Gamma$.~\cite[Proposition 2.2]{CHM2}. 
Among all invariant integer sets $S$ we spot two ones:   
\begin{prop}\label{p.Qmain} For every tile $Q$ and a mask $\bc$, each of the following sets is invariant:

\textbf{a)}  $\ S_0 \ = \ \bigl\{a \in \z^n  \ \bigl| \ |\, (a+Q) \, \cap \, \Gamma \,| \ > \ 0
\bigr\}$; 

\textbf{b)}  $\ \overline S_0 \ = \ \bigl\{a \in \z^n \ \bigl| \  (a+Q) \, \cap \, \Gamma \, \ne \, \emptyset \bigr\}$.  
\end{prop}
{\tt Proof.} We establish a), the proof of b) is the same, replacing 
positivity of the measure by the nonemptiness. 

Take arbitrary $s \in S_0$ and show that if there exist 
$\delta \in \Delta$ and $\gamma \in {\rm supp}\, \bc$ such that the point 
$x = M^{-1} (s + \gamma - \delta)$
is integer, then $x \in S_0$. This will imply that $S_0$ is invariant. 
We have $Mx - s + \delta \, = \, \gamma$. Since the sets 
$s+Q$ and $\Gamma$ possess an intersection of a positive measure, so do the 
shifted sets $s+Q + (Mx - s + \delta)$ and $\gamma + \Gamma$ because 
they are shifted by the same vector.  
Thus, the sets $Q + \delta + Mx$ and $\gamma+ \Gamma$ have an intersection
of positive measure. Hence, so do $M^{-1}(\delta + Q) + x$ and $M^{-1}(\gamma + \Gamma)$. 
However,  $M^{-1}(\delta + Q) \subset G$ and $M^{-1}(\gamma + \Gamma) \, \subset \, \Gamma$. Hence, the bigger sets $x+ Q$ and $\Gamma$ also have an intersection of positive measure, therefore, $x \in S_0$.

{\hfill $\Box$}
\smallskip

For a given finite set $K \subset \z^n$ we denote by $S_K$ the smallest 
invariant set of integers containing~$K$. This  set  is merely 
an intersection of all invariant sets containing~$K$. 

\begin{prop}\label{p.Q0} For the set $K = \{0\}$, 
we have $S_K = S_0$, where the set $S_0$ is defined in Proposition~\ref{p.Qmain}. 
\end{prop}
{\tt Proof.} By Proposition~\ref{p.Qmain}, $S_0$ is an invariant set. 
Clearly, both $Q$ and $\Gamma$ contain a neighbourhood of zero since 
$0\in D$  and $0\in {\rm supp}\, \bc$. Therefore, $S_0$ contains zero and hence  
contains the minimal invariant set $S_K$. Thus, 
$S_K \subset S_0$. If this inclusion is strict, then  
 the set $\Gamma' = Q + S_K$ does not cover $\Gamma$. 
On the other hand, $M^{-1}(\gamma+\Gamma') \, \subset \, \Gamma'$
for all $\gamma \in {\rm supp}\, \bc$. Indeed, for each $q \in S_K$, the set 
 $M^{-1}(q + \gamma + Q)$ is a parallel shift of some $M^{-1}(\delta+Q)$ 
 to an integer vector $x$. Hence $q + \gamma = 
\delta + Mx$ and $x \, \in \, \eta\,  S_K$. 
 Consequently, $x \in S_K$ and $M^{-1}(q + \gamma +Q) \, = \, x\, + \, M^{-1}(d+Q)\, 
 \subset \, x+ \Gamma \, \subset \, \Gamma'$. Thus, $M^{-1}(\gamma + \Gamma') \, \subset \, \Gamma'$. 
 Therefore, the fractal corresponding to the family of contractions $M^{-1}(\gamma \, + \, \cdot \, ), \, \gamma \in {\rm supp}\, \bc\, $ is contained in $\Gamma'$. 
 On the other hand, this fractal is $\Gamma$, and so $\Gamma \subset \Gamma'$, which is 
 a contradiction.   

{\hfill $\Box$}
\smallskip

Proposition~\ref{p.Q0} makes it possible to obtain the set~$S_0$
algorithmically within finite time, without computing the sets~$Q$ and~$\Gamma$. 
This was done in~\cite{CM}, we just slightly modify that construction.  
\bigskip 

\textbf{An algorithm to compute $S_K$ for a given set 
$K$}. 

\noindent {\tt Initialisation.} We have a finite set $K \in \z^n$. Denote $K_0 = K$
and $S_0 = K$. 
The set of digits  $\Delta$ and a mask  $\bc$ are given. 
\smallskip 

\noindent {\tt Main loop.} After the $(j-1)$st iteration we have a finite set of integers $S_{j-1}$ and its subset $K_{j-1}$. Set $\, S_j = S_{j-1}, \, K_j = \emptyset$. For each points 
$k \in K_{j-1}, \,  c \in  {\rm supp} \, \bc$ and $\, \delta \in  \Delta$, we check 
whether or not the point $x = M^{-1}\bigl(k \, + \, c \, - \, \delta\bigr)$ is integer and does not belong to $S_{j-1}$. If so, we set $K_j\, = \, K_j \cup \, 
\{x\}$, otherwise we set 
$K_j\, = \, K_j$. After all triples $(k, c, \delta)$ are exhausted, we set 
$S_j = S_{j-1} \cup \, K_j$. If $K_{j} = \emptyset$, then STOP, 
the algorithm terminates and $S_K = S_j$. Otherwise go to the next iteration. 

\begin{prop}\label{p.QM} For every finite set $K\subset \z^n$, the algorithm 
terminates within finite time and the final set $S_j$ is 
equal to $S_K$.  
\end{prop}
{\tt Proof.} By the construction, $S_j = \cup_{s=0}^j\eta^s\, K$. 
Consider the operator $\xi$ that associates to 
each finite set $K \subset \z^n$ the set  
$\xi \, K \, = \, M^{-1}\bigl(K + {\rm supp} \, \bc \, - \, \Delta\bigr)$. 
Clearly, $\eta K \, \subset \xi K$. Therefore, $S_j$ is contained in the set 
$\cup_{s=0}^{\infty}\xi^s\, K$, whose closure is a fractal set 
of the finitely many contractions $K \, \mapsto \, M^{-1}\bigl(K + c  - \delta\bigr)$, 
where $c \in {\rm supp}\, \bc, \, \delta \, \in \Delta$. Since this fractal set is compact, it contains only a finite number of integers. All the sets $S_j$ produced by the algorithm are  contained in this finite set. Hence, for some $j$ we necessarily have $S_j = S_{j-1}$. Since $S_j = \cup_{s=0}^j\eta^s\, K \, = \, 
\cup_{s=0}^{\infty}\eta^s\, K$, we see that $\eta\, S_j \, \subset S_j$, 
so $S_j$ is invariant. On the other hand, each invariant set that contains 
$K$ must also contain $\eta^s  K$  for all $s$, hence it contains $S_j$. 
Thus $S_j = S_K$.

{\hfill $\Box$}
\smallskip

\bigskip 

\begin{center}
\textbf{4.4 Spectral factorisation of  the dilation matrix}
\end{center}
\bigskip 

Now we need to use spectral 
properties of the dilation matrix~$M$. 
All eigenvalues of~$M$ are bigger than one by modulus. 
Let $r_1 < \cdots < r_q$ be all possible absolute values 
of eigenvalues of $M$ and let exactly $n_i$ of them (counting multiplicities) be equal 
by modulus to~$r_i, \, i = 1, \ldots , q(M)$.  Let $J_i\subset \re^n$ be the linear span of root
subspaces of~$M$ corresponding to all eigenvalues of modulus~$r_i$. Thus, $\hbox{dim}(J_i)=n_i$ and
the operator~$M|_{J_i}$ has all its eigenvalues equal to $r_i$ in the absolute value.
The space $\re^n$ is a direct sum of $J_1, \ldots , J_q$: 
$$
 \re^n\quad =\quad \bigoplus_{i=1}^{q} \, J_i\, . 
$$
There exists an invertible transformation $B: \re^n \rightarrow \re^n$ such
that $M$ has the following block diagonal structure
\begin{equation}\label{eq.M-factor}
B^{-1} M B\quad = \quad
\left(
\begin{array}{cccc}
 M|_{J_1} & 0 & \cdots & 0\\
 0       & M|_{J_2} &      & \vdots \\
\vdots   &         & \ddots & 0 \\
0 & \cdots  & 0 & M|_{J_q}
\end{array}
\right).
\end{equation}
The subspaces~$J_k$ will be referred to as {\em spectral subspaces}
and~(\ref{eq.M-factor}) is a {\em spectral factorisation}. 
In particular, if the matrix~$M$ is isotropic, then 
$q(M)=1$, in which case $J_1$ coincides with~$\re^n$. 
The converse is not true: if all eigenvalues of $M$ have equal moduli and 
the~$M$ has nontrivial Jordan blocks, then it is not isotropic, although~$q=1$. 
\bigskip 

\begin{center}
\textbf{4.5 Wide simplex and admissible sets}
\end{center}
\bigskip 

The invariant set  $S_0$ defined in Proposition~\ref{p.Qmain} 
possesses the following key property: {\em shifts of the tile~$Q$ over vectors from 
$S_0$ cover the set~$\Gamma$.} The set $\overline S_0$ has more: {\em shifts of $Q$ 
over  vectors from  $\overline S_0$
cover a neighbourhood of~$\Gamma$.} We need to define a property which is between those two. First we introduce two more notation. 

A {\em wide simplex} is a simplex in~$\re^n$ with one of vertices 
at the origin such that its interior intersects all spectral subspaces~$J_k, 
k = 1, \ldots , q$, 
of the matrix~$M$. 

The existence of wide simplices is easily shown (see also~\cite{CP}). 
Moreover, a homothety about the origin respects wide simplices. Hence, 
every ball centered at the origin contains a wide simplex. 
Actually, even every half-ball contains a wide simplex.   

\begin{defi}\label{d.admiss}
Let $\Gamma$ be a support set of refinement equation defined by~(\ref{eq.K})
and let $Q = Q(M, \Delta)$ be a tile.  A finite subset~$S \subset \z^n$ is called admissible if the set $S+Q$ contains the sum of $\Gamma$ 
with some wide simplex. 
\end{defi}
Since every ball centered at the origin contains a wide simplex, 
we see that $S$ is admissible whenever $S+Q$ contains a neighbourhood of $\Gamma$. 
In particular, the set~$\overline S_0$  from Proposition~\ref{p.Qmain} is 
always admissible. The set $S_0$ may be not, however, in most cases it is  admissible
as well. Actually we can always use the set~$S = \overline S_0$. However, in some cases it is too large and can be replaced by a smaller admissible set. 
This is very important from the computational point of view, see~Remark~\ref{r.35}.
This is the only reason for introducing the notions of wide simplices and 
of admissible sets.   

\bigskip 

\begin{center}
\textbf{4.6 The vector-function $v(x)$ and the transition matrices 
$T_{\delta}, \, \delta \in \Delta$}
\end{center}
\bigskip 

Now we are realising the main idea of the matrix approach to 
multivariate refinement equations. For an arbitrary refinement equation~(\ref{eq.ref}), we first  choose a tile $\, Q = Q(M, \Delta)$
(basic tile) and then 
we pass from the function~$\varphi: \, \re^n \, \to \, \re$ 
to the vector-valued function~$v: Q \, \to \, \re^N$. To define 
this function, we take an arbitrary admissible invariant set $S \subset \z^n$
(Definition~\ref{d.admiss}) and denote $|S| = N$. Then 
$v(x)$ is defined as follows:  

\begin{equation}\label{eq.v}
v\, :\, Q \rightarrow \re^N, \quad v(x) \ = \ v_{\varphi}(x) \ = \ 
\Bigl(\, \varphi(x+k)\, \Bigr)_{k \in S}\, , \quad x \in Q\, .
\end{equation}
For convenience, we enumerate the components of the vector~$v$ by elements of 
the set~$S$. Consider the $m$ following $N\times N$ {\em transition matrices}
$T_{\delta}, \, \delta \in \Delta$,  defined by the equality 
\begin{equation}\label{eq.T}
(T_{\delta})_{a b} \ = \ c_{\, Ma - b + \delta}\ , \quad a, b \, \in \, S\, , \quad 
\delta \in \Delta.
\end{equation}
Rows and columns of the transition matrices are enumerated by elements of  the set~$S$. 
We denote $\, \cT \, = \, \{T_{\delta} \ : \ d \in \Delta\}.$
The refinement equation on the function~$\varphi(x)$ is equivalent to
   the following equation on 
the vector-valued function~$v(x)$: 
\begin{equation}\label{eq.ss1}
v (x) \ = \ T_{\delta} \, v(Mx - \delta)\ , \qquad x \in M^{-1}(Q + \delta)\ , \quad \delta \in \Delta\, .
\end{equation}
Functional equations of this type are often called {\em equation of self-similarity}~\cite{P08}. 
\begin{remark}\label{r.25}
{\em If in the definition of $v(x)$
we used an arbitrary invariant set~$S$, then (rather surprisingly!) the   
$L_p$-regularity of~$v$ might 
not be equal  to the $L_p$-regularity 
of~$\varphi$.  The reason is that the function~$v$ is defined 
on~$Q$, while $\varphi$ is defined on the entire~$\re^n$. That is why 
we had to use an admissible invariant set~$S$. This guarantees that 
the union of translations $\bigcup_{k \in S}(k+ Q)$ contains not only the support of $\varphi$ but a bigger set: the sum of this support with a wide simplex. Let us note that 
for measuring the H\"older regularity in~$C(\re^n)$, this enlargement is not needed
and  any invariant set~$S$ 
suffices~\cite{CP}.  The reason is that a continuous refinable function vanishes 
on the boundary of its support, which may not be true for an $L_p$ refinable function.   }
\end{remark}

\bigskip

\begin{center}
\textbf{4.7 Special subspaces of $\re^N$}
\end{center}
\bigskip 

In the regularity analysis of refinable functions we deal with several linear 
affine subspaces in~$\re^N$. First we define 
$$
 V\ =\ \Bigl\{\, w\, =\, (w_1, \ldots, w_N) \in \re^N \quad : \quad \sum_{j=1}^{N} w_j\, =\, 1\, \Bigr\}.
$$
It is well known that every compactly supported refinable function
such that $\int_{\re^n} \varphi(x)\, ds = 1$ possesses the 
{\em partition of unity property:}
$$
\sum_{k \in \z^n} \varphi (x+k)\quad  \equiv \quad  1
$$
Hence, after a multiplication of $\varphi$ by a constant it may be assumed that $v(x) \in V$ for almost all $x \in G$.  
We denote  the linear part of the affine subspace~$V$ by  
$$
W\quad =\quad \Bigl\{\, w\ =\ (w_1, \ldots, w_N) \in \re^N \quad : \quad  \sum_{j=1}^N w_j=0 \, 
\bigr\}\, . 
$$
Finally,  define the space of differences of the vector-function 
$v = v_{\varphi}$: 
\begin{equation}\label{eq.U}
 U\ = \ {\rm span}\, \Bigl\{\, v(y)\, - \, v(x) \quad :  \quad y, x \, \in Q\, \Bigr\}\, . 
\end{equation}
Since $v(x) \in V$ for almost all $x \in Q$, we  have  $U \subset W$. 
The sum rules~(\ref{eq.sumrules}) imply that the column sums of each matrix $T_{\delta}$ are equal to one. 
Therefore, $T_{\delta}V \subset V$ and $T_{\delta} W \subset W$
for all~$\delta \in \Delta$. Thus, $V$ is a common affine invariant subspace of the family 
$\cT$ and $W$ is its common linear invariant subspace. 

For $i=1, \ldots , q$, define the subspaces $U_1, \ldots U_q$
of the space $\re^N$ as follows: 
\begin{equation}\label{eq.ui}
U_i\ = \
{\rm span}\, \Bigl\{v(y)\, - \, v(x) \quad :  \quad  x, y \in Q, \  y-x \in J_i \Bigr\}\, , \quad
i=1, \ldots, q(M).
\end{equation}
Note that $U_i$ are nonempty, due to the interior of $Q$ being nonempty. It is seen easily that the spaces $\{U_i\}_{i=1}^{q}$ span the whole space $U$,
but their sum may not be direct. The subspaces $\{U_i\}_{i=1}^{q}$, unlike the subspaces $\{J_i\}_{i=1}^{q}$, may have nontrivial intersections. For example, they can all coincide with $U$. It turns out that all $U_i$ are common invariant 
subspaces for the matrices~$T_{\delta}$. 

\begin{lemma}\label{l.invar}
If $J$ is an invariant subspace for the matrix $M$, then
$L \, = \, {\rm span}\, \{v(y) - v(x) \ : \ y-x \in J\}$ is a common invariant subspace for all $T_{\delta}, \, \delta \in \Delta$.
\end{lemma}
{\tt Proof}. If $u \in L$, then $u$ is a linear combination 
of several vectors of the form $v(y) - v(x)$ with $y - x  \in J$.
For every $\delta \in \Delta$ we define $x' = M^{-1}(x+\delta), y' = M^{-1}(y+\delta)$ and have
$$
v(y') - v(x') \ = \ T_{\delta} \, \bigl( v(My'-\delta) \, - \, v(Mx' - \delta)\bigr) \ = \ T_{\delta} \, \bigl( v(y) \, - \, v(x)\bigr)\, .
$$
Hence, $T_{\delta}\bigl( v(y) -  v(x)\bigr) \in L$ for each pair $(x, y)$, and, therefore, $T_{\delta} u \in L$ for all $u \in L$.

{\hfill $\Box$}

\smallskip

\smallskip

\begin{center}
\textbf{4.8 The formula of regularity for refinable 
functions in $L_p$}
\end{center}
\bigskip 

This formula expresses the H\"older exponent of the refinable function with 
the $L_p$-spectral radius of matrices~$T_{\delta}$ restricted to the subspaces~$U_i$. 
For a given set of linear operators $\cA = \{A_0, \ldots , A_{m-1}\}$ acting in~$\re^d$ and for given $p \in [1, +\infty)$, the $L_p$-spectral radius ($p$-radius) 
is defined by the formula: 
$$
\rho_p\quad =\quad \rho_p(\cA)\ = \ \lim_{k \to \infty}\, \Bigl(\,
    m^{-k}\sum_{A_{\ell_i} \in \cA, \, i = 1, \ldots , k}\|A_{\ell_1}\cdots A_{\ell_k}\|^p \, \Bigr)^{\, 1/pk}.
$$
The limit always exists and does not depend on the operator 
norm (see~\cite{P97}.  
for more on properties of the $p$-radius). Clearly, for one operator, the value~$\rho_p$
becomes the usual spectral radius, i.e., the largest by modulus  eigenvalue.
Already for two operators, the computation of the $p$-radius is a hard problem.  
For example, it is still not clear if the $1$-radius can be efficiently computed. 
On the other hand, for even integer~$p$, the $p$-radius can be expressed 
by means of a usual spectral radius of some large matrix. For example, 
the $2$-radius is equal to the square root of the spectral radius of 
the following operator~$\mathbf{A}$ acting on the space~$\cM_d$ of symmetric $d\times d$-matrices: 
\begin{equation}\label{eq.A}
\mathbf{A}(X) \ = \ \frac{1}{m}\, \sum_{i=0}^{m-1} A^*_iXA_i\, , \qquad X \in \cM_d\, . 
\end{equation}
This operator acts on the $\frac{d(d+1)}{2}$-dimensional space $\cM_d$
and obeys an invariant cone of positive semidefinite matrices. Hence, 
by the Krein-Rutman theorem~\cite{KR}, its largest by modulus eigenvalue~$\lambda_{\max}$ (which can also be called Perron eigenvalue)
is positive. 
The fact is $\rho_2(A_0, \ldots , A_{m-1})\, = \, \sqrt{\lambda_{\max}}$~\cite{LauW, P97}.

The formula for $L_p$-regularity of univariate refinable functions 
was well-known~\cite{LW, P97}. However, it offered a surprising resistance in extending to multivariate functions. For general dilation matrices~$M$, this 
extension was done only in 2019~\cite{CP}. The main idea is to find  the H\"older exponent
separately on the spectral subspaces~$J_i$. 
The H\"older exponent of $\varphi$ along a subspace $J \subset \re^n$ is defined by
$$
\alpha_{p, J}(\varphi)\quad = \quad \sup
\Bigl\{\, \alpha \ge 0 \quad : \quad \bigl\|\varphi(\cdot  + h) \, - \, \varphi(\cdot)\bigr\|_p \ \le \  C\,  \|h\|^{\, \alpha}\, , \quad h \in J\, \Bigr\}\,.
$$
The following theorem was proved in~\cite{CP}. Let us remember that 
$\cT  = \{T_{\delta} \, , \, d \in \Delta\}.$
\begin{theorem}\label{th.holder-p} Let $1 \le p <\infty$.
For a refinable function~$\varphi \in L_p(\re^n)$, we have
\begin{equation}\label{eq.holder-direct-p}
\bal_{p, J_i}(\varphi)\quad = \quad  \log_{\, 1/r_i} \, \rho_{p}(\cT|_{{U}_i})\, , \qquad  i = 1, \ldots , q\,
\end{equation}
and, consequently,
\begin{equation}\label{eq.holder-p}
\bal_{p}(\varphi)\quad = \quad   \min\limits_{i = 1, \ldots , q}\,
\log_{\, 1/r_i} \, \rho_{p}(\cT|_{{U}_i})
\end{equation}
\end{theorem}
For isotropic matrices, when all $r_i$ are equal to $r = \rho(M)$, 
formula~(\ref{eq.holder-p}) looks as simple as for the univariate 
refinable functions: $\bal_{p}(\varphi)\, = \,  
\log_{\, 1/r} \, \rho_{p}(\cT|_{{U}})$. 

\begin{remark}\label{r.35}{\em Theorem~\ref{th.holder-p} 
expresses the $L_p$-H\"older regularity of a refinable function 
to $p$-radii of the transition matrices~$T_{\delta}$ restricted to 
special common invariant subspaces. As we have mentioned, for even integer~$p$, 
the $p$-radius can be computed as a Perron eigenvalue of some high-dimensional matrix, for other $p$ only approximate computational methods are known. 
At any rate, the complexity of computation depends significantly of the size of matrices $T_{\delta}$, which is $N = |S|$. That is why it is important to 
reduce the cardinality of the admissible set~$S$. The set 
$S = \overline S_0$ is sometimes too large and it is possible to find a smaller 
admissible set using the notion of wide simplices.    

}
\end{remark}

\newpage

\bigskip 

\begin{center}
\textbf{5. Computing the surface regularity \\ and surface dimension of attractors and tiles}
\end{center}
\bigskip 

The characteristic function of an attractor satisfies functional equation~(\ref{eq.ref0}), which is the refinement equation with 
the coefficients  $c_k = 1$ if $k\in D$ and $c_k = 0$ otherwise. 
Therefore, Theorem~\ref{eq.holder-direct-p} can be applied directly for 
computing the H\"older  regularity and (if the dilation matrix~$M$ is isotropic)
the surface regularity and the surface dimension of attractors. 
The specific mask containing only zeros and ones makes the computation 
easier. Moreover, it will allow us to come up with simpler 
formulas of regularity that do not involve the subspaces~$U_i$, 
which are a priori, unknown (Theorem~\ref{th.tile-isotr}). 
This means that the same formulas can be applied to find 
the H\"older exponents and the surface regularity of multivariate Haar wavelets 
generated by arbitrary dilation matrices.

First of all, we observe that the refinement equation
for attractors~(\ref{eq.ref0})  admits 
the transition matrices~$T_{\delta}$, which 
will be {\em simple} in the following sense:   

\begin{defi}\label{d.simple}
A matrix is called simple if each its column contains precisely one entry equal to one and all others are zeros.
\end{defi}

\begin{prop}\label{p.simple}
Suppose $G(M, D) \subset \re^n$ is an attractor;  
then for every basic tile $Q(M, \Delta)$ and for every admissible invariant set $S$, 
the matrices $T_{\delta}, \, \delta \in \Delta$,  are all simple. 
\end{prop}
{\tt Proof}. We have~$c_k = 1$ if and only if $k \in D$, otherwise 
$c_k = 0$. Therefore (formula~(\ref{eq.T})), 
$(T_{\delta})_{a b} =  c_{\, Ma - b + \delta} = 1$ 
if and only if   $Ma - b + \delta \, \in \, D$. 
If $b$ is fixed, then the set $\, b \, - \, \delta \, + \, D$
has precisely one common point with the 
lattice $M\, \z^n$, since $D$ is a digit set. Therefore, one 
component of the $b$th column of $T_{\delta}$ is one and the others are zeros. 

{\hfill $\Box$}
\smallskip 

All simple matrices form  a miltiplicative matrix semigroup. 
Let us remember that $V$ is an affine hyperspace of~$\re^N$ 
that consists of points with the sum of components  being one and 
$W$ is its linear part. 
A simple matrix is column-stochastic, hence it respects both $V$ and $W$.
As usual,  $\cA^k$  denotes the set of all  products of matrices from~$\cA$
of length~$k\ge \n$ (repetitions permitted). Clearly, $|\cA^k| = m^k$. We write $\cA^k_0$ for the set of matrices from~$\cA^k$ that have at most one positive 
row. If all matrices from $\cA$  are simple, then so are all matrices 
from~$\cA^k$ and each matrix from~$\cA^k_0$ has one row of ones and all other elements are zeros.   
\begin{prop}\label{p.s}
Let $\cA =  \{A_0, \ldots , A_{m-1}\}$ be a set of simple matrices. 
Then for every common invariant subspace $U \subset W$ of the matrices
from $\cA$,  we have 
$\rho_p \, (\cA|_{U}) \ = \  \bigl[ \rho_1 \, (\cA|_{U})\bigr]^{1/p}$. 
In case $U = W$, the following formula  holds: 
\begin{equation}\label{eq.row}
\rho_1 \, (\cA|_{U}) \ = \ \lim_{k\to \infty }\ \left[\, 1 \ - \ 
\frac{|\cA^k_0|}{|\cA^k|}\, \right]^{\, 1/k} 
\end{equation}
\end{prop}
{\tt Proof}. We begin with the case $U=W$. 
The norm of every simple matrix restricted to $W$
is either zero (if this matrix has precisely one non-zero row) or 
between $1$ and $\sqrt{\frac{n}{2}}$  otherwise. 
Replacing the norms of all matrix products in 
the definition of $L_1$-spectral radius by those numbers we see that 
the quantity~(\ref{eq.row}) is equal to $\rho_1(\cA|_{W})$. 
Now consider the case of general $U \subset W$. Denote  by $\cH$ the set 
of all simple $N\times N$ matrices. 
This set is finite and for every~$k$, all product from $\cA^k$
belong to this set. Hence, the number  $\|\Pi|_{U}\|$
for arbitrary $\Pi \in \cA^k$ and $k \in \n$, 
can take a finite number of values. Consequently, 
 $\|\Pi|_{U}\|^p \, \asymp \, \|\Pi|_{U}\|$
 and this equivalence is defined by two absolute constants. 
 Therefore, the values $\, m^{-k}\sum_{\Pi \in \cA^k} \|\Pi|_{U}\|^p\, $ and 
 $\, m^{-k}\sum_{\Pi \in \cA^k} \|\Pi|_{U}\|^p\, $ are equivalent by the same two constants. Hence, $\rho_p \, (\cA|_{U}) \ = \  \bigl[ \rho_1 \, (\cA|_{U})\bigr]^{1/p}$, which concludes the proof.

{\hfill $\Box$}
\smallskip 

\bigskip

Applying Theorem~\ref{th.holder-p} and Proposition~\ref{p.s} we obtain
\begin{theorem}\label{th.attr}
For every attractor~$G$, we have 
$$
 \bal \ = \   2\, \min_{1, \ldots , q}\, 
 \log_{\, \frac{1}{r_i}}\,  \rho_2\,(\cT|_{U_i})\, .
$$ 
If the matrix $M$ is isotropic, then  
$$
\bs \ = \ \bal \ = \   2\, \log_{\, \frac{1}{r}}\,  \rho_2\,(\cT|_{U})\, .
$$ 
\end{theorem}
This theorem allows us to find the H\"older regularity of any attractor 
as the largest eigenvalue of the operator~$\mathbf{A}$ defined by~(\ref{eq.row})
for $m$~operators~$A_{\delta} \, = \, T_{\delta}|_{U}, \, \delta \in \Delta$. 
If in addition the dilation matrix is isotropic, then the surface regularity
is equal to the same value. In case $G$ is a tile, everything can be rewritten in a simpler terms. Moreover, in this case we can use the fact that all matrices $T_{\delta}$ are simple and the $L_1$-spectral radius~$\rho_1\,(\cT|_{W})$ obtains a combinatorial form~(\ref{eq.row}).

\begin{theorem}\label{th.tile-isotr}
For every time $G$  generated by an isotropic matrix~$M$, 
we have 
$$
\bs \ = \ \bal \ = \   2\,\log_{\, \frac{1}{r}}\,  \rho_2\,(\cT|_{W})\, .
$$ 
This value is equal to the right hand side of equality~(\ref{eq.row}). 
\end{theorem}
The main advantage of this theorem is that the expression for the 
H\"older exponent and for the surface regularity of a tile does not depend on the 
subspace~$U$ (which may be different for different tiles) and can be expressed 
by the $1$-radius of the matrices~$T_{\delta}$ restricted to the standard 
subspace~$W \, = \, \{w \in \re^N \, | \,  \sum_{i}w_i = 0\}$. The proof requires one auxiliary result. As usual we denote by $\cT^k$ the set of all 
products of length~$k$ of matrices from~$\cT$.

\begin{lemma}\label{l.20}
 Let $G(M, D)$ be an attractor and $S\subset \z^d$ be an admissible set 
 for $G$; then there is $p\in N$ such that 
 for every $i\in S \setminus \overline S_0$, all matrices from~$\cT^p$
 have zero $i$th row. 
\end{lemma}
{\tt Proof}. Consider the {\em transition operator}
$$
\bT\,f(x)\ = \ \sum_{k \in D} f(Mx - k)\, .
$$ 
Clearly, for the function~$\varphi = \chi_{G}$, we have $\bT\varphi = \varphi$. 
If $f$ is the indicator function of some compact set $K \subset \re^n$, 
then the Hausdorff distance between the support of the function~$\bT^pf$
and $G$ tends to zero as $p\to \infty$. Indeed, since $M^{-1}$ has spectral radius 
smaller than one, there exists a norm in $\re^n$ such that in the corresponding 
operator norm~$q = \|M^{-1}\| < 1$. In this norm, 
the distance between the support of~$\bT^pf$ and $G$ is at most 
$q^j$ times the distance between the support of~$f$ and $G$.  
Since $i\notin \overline S_0$, it follows that $i + Q$ does not intersect~$G$. 
Therefore, for all sufficiently large $p$, the support of~$\bT^pf$
does not intersect the set $i + Q$. Let now $K = j+Q$ for some $j\in S$ and 
$f = \chi_{K}$.  
Then   for every sequence $\ell_1, \ldots , \ell_k$, 
the element in the $i$th row and $j$th column of the 
matrix $\Pi = T_{\ell_1}\cdots T_{\ell_p}$ is equal to 
$\bT^pf(i+0.\ell_1\ldots \ell_k)$. However, $i+0.d_1\ldots d_k\, \in \, 
i + Q$, hence $\bT^pf(i+0.\ell_1\ldots \ell_k)$ and so $\Pi_{ij} = 0$. 
Thus, for all sufficiently large~$i$, the $i$th row of every product~$\Pi \in \cT^p$ is zero. 
 
{\hfill $\Box$}
\smallskip 

{\tt Proof of Theorem~\ref{th.tile-isotr}}. 
Since $G$ is a tile,  for almost all $x$, the vector $v(x)$ has only one non-zero component, which is equal to $1$. Namely, $v_k =1$ if $x \in k+G$. 
Denote $\cI \in S \cap \overline S_0$  
and $L \, = \, {\rm span}\,\{e_i -e_j \ | \ i,j \, \in \, \cI\}$. 
Then $U = L$. Furthermore, by Lemma~\ref{l.20}, there is 
$p$ such that for all~$i \notin \cI$, 
the $i$th row of every matrix from~$\cT^p$ is zero.   Therefore, 
$\rho_1(\cT|_{W}) \, = \, \rho_1^{1/p}(\cT^p|_{W}) \, = \, 
\rho_1^{1/p}(\cT^p|_{L})\, = \, \rho_1(\cT|_{L})$ and hence 
$\rho_1(\cT|_{U}) \, = \, \rho_1(\cT|_{W})$. 
On the other hand, since $M$ is isotropic, then all 
$r_i$ are equal to~$r$, therefore formula~(\ref{eq.holder-p})
reads $\bal \, = \, \log_{\, \frac{1}{r}}\, \rho_1\, (\cT|_{U})$. 
This completes the proof.

{\hfill $\Box$}
\smallskip

\newpage 

\begin{center}
\textbf{6. Examples}
\end{center}

\begin{ex}\label{ex.20}
{\em {\em The univariate tile} $G(M, D)$ with $M=3$ and $D = \{0,1,5\}$. 
Its characteristic function $\varphi = \chi_{G}$ satisfies the refinement equation
$$
\varphi(x) \ = \ \varphi(3x) \ + \ \varphi(3x-1)\ + \ \varphi(3x-5), \quad x \in \re,
$$
Since every point of $G$ has the form 
$x\, =\, \sum_{k=1}^{\infty}\ell_k 3^{-k}$ with $\ell_k \in D$, we have 
${x \, \le \, 5\, \sum_{k=1}^{\infty} 3^{-k}\, \le \, \frac52}$.  
Therefore $G \subset [0, 2.5]$. This tile is shown in fig.~1. 

\begin{figure}

\centering

\includegraphics[scale=0.3, bb=10 -40 1433 430]{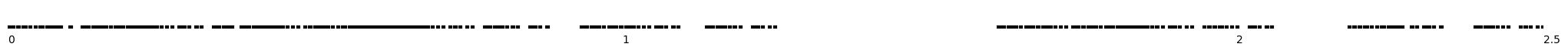}

\caption{ \footnotesize{The one-dimensional tile with 
$M = 3, D=\{0,1,5\}$:  $\bs = 0.1977..$ 
and $\bd = 0.8022..$}}

\label{fig1}

\end{figure}

\smallskip

\noindent Taking the basic tile $Q = [0,1]$ 
generated by the digit set  $\Delta = \{0, 1, 2\}$, we have
$\cup_{k=0,1,2} (k+Q)\, = \, [0,3]$. 
Therefore, $S = \{0,1,2\}$ is an admissible invariant set. 
Hence, $N = 3$ and there are three transition matrices: 
$$
T_0 \ = \ \left(
\begin{array}{ccc}
1 & 0 & 0\\
0 & 0 & 1\\
0 & 1 & 0
\end{array}
\right)\ , \quad
T_1 \ = \ \left(
\begin{array}{ccc}
1 & 1 & 0\\
0 & 0 & 0\\
0 & 0 & 1
\end{array}
\right)\  \quad \hbox{and} \quad
T_2 \ = \ \left(
\begin{array}{ccc}
0 & 1 & 1\\
1 & 0 & 0\\
0 & 0 & 0
\end{array}
\right).
$$
The subspace $W$ is two-dimensional. Choosing the 
basis $e_1 = (1, -1, 0)^T, e_2 = (0, 1, -1)^T$ of $W$, 
we obtain  the matrices $A_{\delta}=T_{\delta}|_W$, $\delta = 0,1,2 $: 
$$
A_0 \ = \ \left(
\begin{array}{cc}
1 & 0 \\
1 & -1
\end{array}
\right)\ , \quad
A_1 \ = \ \left(
\begin{array}{rr}
0 & 1\\
0 & 1
\end{array}
\right)\  \quad \hbox{and} \quad
A_2 \ = \ \left(
\begin{array}{rr}
-1 & 0 \\
0 & 0
\end{array}
\right)
$$
To compute $\rho(\cT|_W) = \rho_2 (A_0, A_1, A_2)$ we build the matrix of the 
operator  $\mathbf{A}$ defined by~(\ref{eq.A}). The 
space of symmetric $2\times 2$ matrices has dimension~$3$, and 
$$
\mathbf{A} \ = \ 
\frac13\, 
\left(
\begin{array}{ccc}
2 & 1 & 2\\
1 & 2 & 2\\
0 & -1 & -1
\end{array}
\right).
$$
For this matrix, $\lambda_{\max} \, = \, \frac{1+\sqrt{2}}{3}$, 
and therefore~$\rho_2(\cT|_W) \, = \, \sqrt{\frac{1+\sqrt{2}}{3}}$. 
Now by Theorem~\ref{th.tile-isotr}, we have $\bs(G) = \bal(G) = 
-\log_3 \frac{1+\sqrt{2}}{3} \, = \, 0.1977...$.  Hence, the surface dimension 
of~$G$ is $\, \bd \, = \, 1 - \bs \, = \, 0.8022...$. 
}
\end{ex}

\begin{figure}

\centering

\includegraphics[scale=0.25, bb=5 10 733 430]{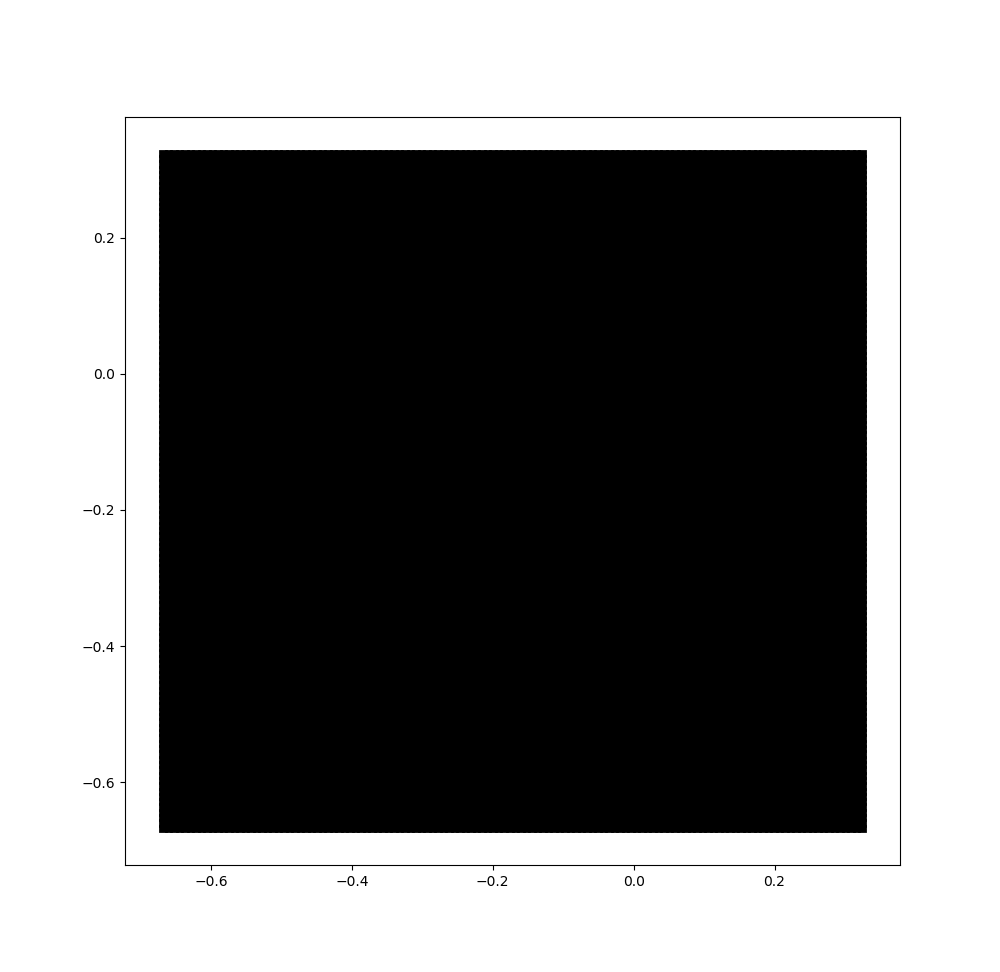}

\caption{\footnotesize{The Square:  $\bs = 1$ 
and $\bd = 1$}}
\vspace{2cm}

\end{figure}

\begin{ex}\label{ex.30} Two-digit tiles on the plane.  
{\em There are only three attractors in~$\re^2$ up to affine similarity
which are generated by two digits, i.e., when $m=2$, see~\cite{Z1}. 
In all the three cases the digit set can be~$D=\{(0,0)\, ; \, (1,0)\}$. 
In this case all those three attractors are tiles. 

The first one 
is a unit square (fig.~2), it is generated by the matrix~$M\, = \, \left(0 \, -2 \atop 1\ \,  0 \right)$. Of course, for this tile $\bs = 1$ and $\bd = 1$. 

\smallskip

The second type is more interesting. This is the Dragon (fig.~3) generated 
by the matrix 
$$
M \ = \ \left(
\begin{array}{rr}
\, 1 & 1 \\
-1 & 1
\end{array}
\right)
$$
(rotation on~$45^o$ with the expanding by~$\sqrt{2}$). Denote this tile by~$G$.

To compute $\bs$ we choose the set $\overline S_0$ (Proposition~\ref{p.Qmain}), 
in which case the basis tile~$Q$ coincides with~$G$. We have 
$\overline S_0 \, = \, (0,0), \, (\pm 1, 0)\, , \, (0, \pm 1)\, , \, (\pm 1, \pm 1)$
(seven points), therefore $N=7$ and the matrices $T_0$ and $T_1$ are $7\times 7$.
Computing~$\rho_2(\cT|_W)$ we obtain $\bal = 0.4763..$.  Therefore $\bs = 0.4763.. $
and $\bd = 1.5236..$.  

\begin{figure}

\centering

\includegraphics[scale=0.25, bb=5 10 733 430]{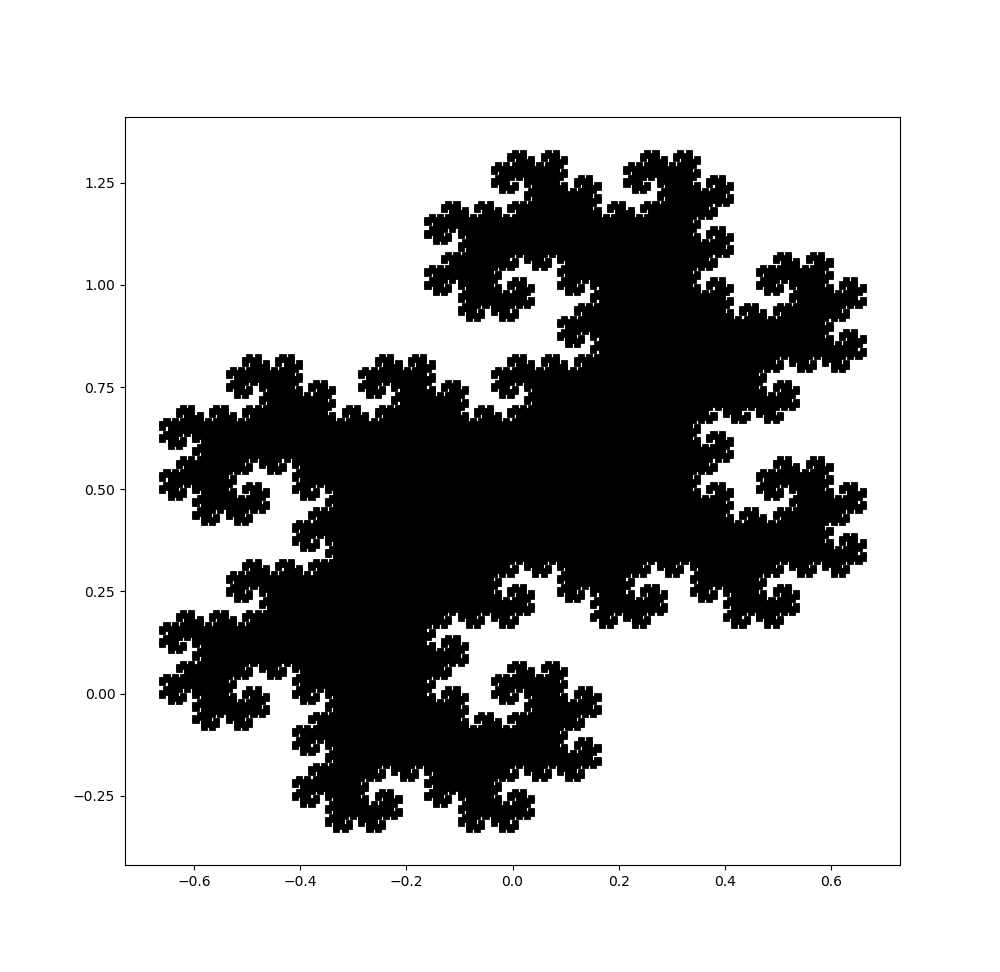}

\caption{\footnotesize{The Dragon:  $\bs = 0.4763.. $ 
and $\bd = 1.5236..$}}

\end{figure}

\newpage 

\smallskip 

The third two-digit tile is generated by the matrix 
$$
M \ = \ \left(
\begin{array}{rr}
\, 1 & -2 \\
1 & \, 0
\end{array}
\right)
$$
\begin{figure}

\centering

\includegraphics[scale=0.3, bb=5 10 733 430]{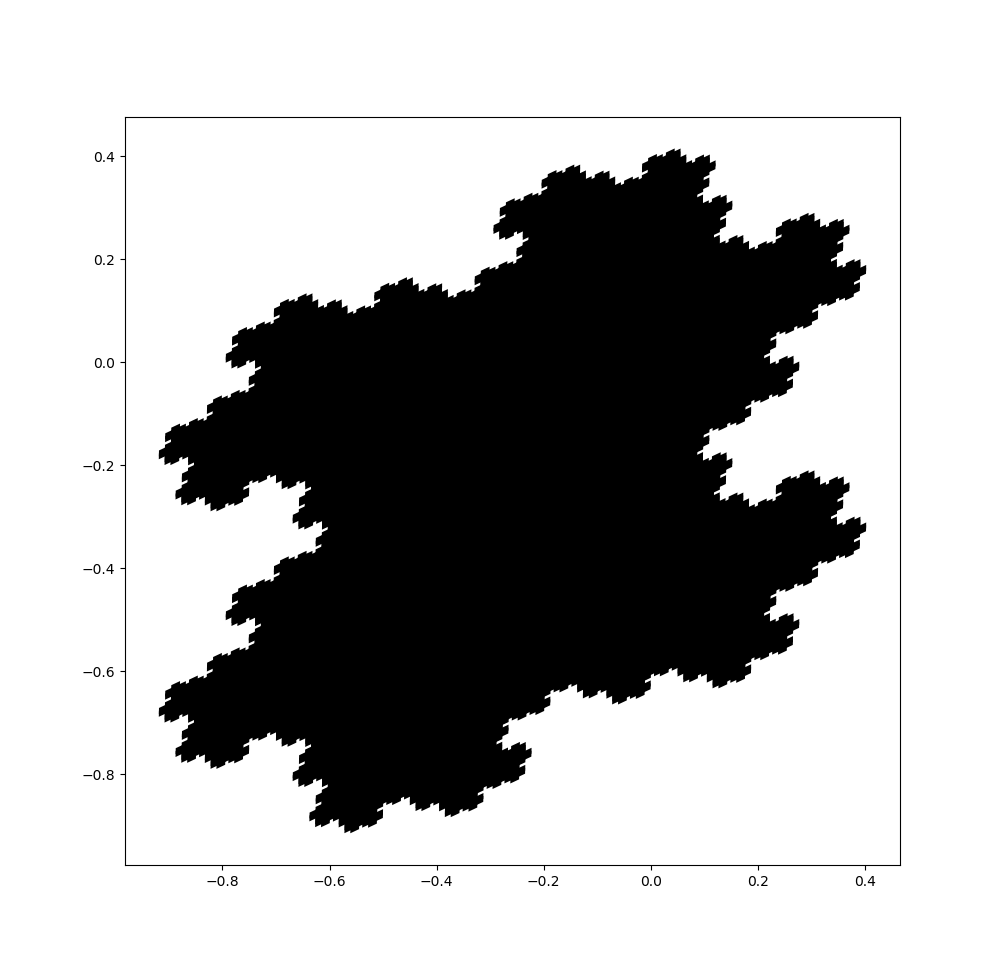}

\caption{\footnotesize{The Bear:  $\bs = 0.7892.. $ 
and $\bd = 1.2107..$}}
\vspace{2cm}


\end{figure}

\smallskip 
As it is seen in fig.~4, if is natural to call it {\em Bear}. 
It has the same set $\overline S_0$ and hence its transition matrices~$T_0, T_1$
are again $7\times 7$.  Computing~$\rho_2(\cT|_W)$ we obtain $\bal =0.7892..$.  Therefore $\bs = 0.7892..$
and $\bd = 1.2107..$. Thus, the Bear has a bigger surface regularity than 
the Dragon. 

}
\end{ex}

\begin{ex}\label{ex.40} Plane tiles with $M = 2I$.  
{\em In this case
$\, m = |\det M|= 4$, hence there will be four digits. 
Different choices of these digits define different attractors. 
For example, for $D\, = \, \bigl\{(0,0)\, ; \,  (1,0)\, ; \, 
(1,0)\, ; \, (1,1) \bigr\}$, we obtain a unit square. 
In this case, of course, $\bal = 1$ and $\bd = 1$. 
Changing one digit: $(1,1)$ to $(-1, -1)$, we obtain the tile 
depicted in fig.~5

  For this tile $\bs \, = \, 0.4150...$
and respectively $\bd \, = \, 1.5849...$. So, its regularity is close to 
the Dragon. It looks similar to the Sierpinski 
carpet ans can be called {\em quasi Sierpinski tile}. In contrast to the 
Sierpinski carpet, which has measure zero, it has measure~$1$ as a tile. 
This tile was considered in~\cite{W}. 
}
\end{ex}

\begin{figure}

\centering

\includegraphics[scale=0.4, bb=5 10 733 430]{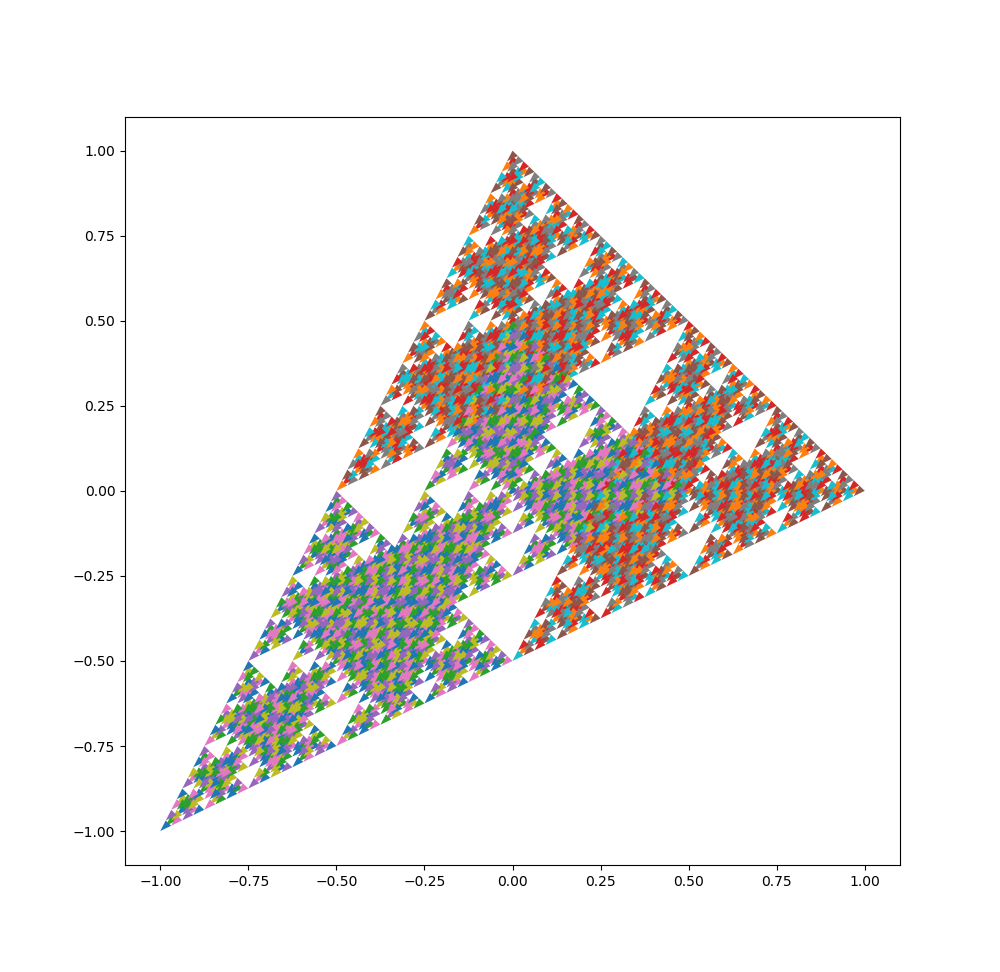}

\caption{\footnotesize{The quasi-Sierpinski tile. $\bs \, = \, 0.4150..., \, 
\bd \, = \, 1.5849...$}}


\end{figure}

\begin{center}
\textbf{7. Attractors of the highest regularity}
\end{center}
\bigskip 

The highest possible surface regularity of any set is one. It is attained, for example, 
 for sets bounded by surfaces of finite area: for convex sets,  
for sets 
 with piecewise smooth boundaries, or for finite unions of such sets.  
 We conjecture that for self-affine attractors, this situation is impossible 
 apart from the case of parallelepipeds.    

\begin{conj}\label{conj.10}
If an attractor satisfies $\bs = 1$, then it is a parallelepiped. 
\end{conj}

For two-digit attractors on the plane this is true, since there are only three types of such attractors~\cite{Z1}. In Example~\ref{ex.30} we analysed all of them: 
for the square we have~$\bs = 1$, for Dragon and for Bear, $\bs < 1$. In every dimension~$n$, for each $m\ge 2$, there are finitely many, up to an affine similarity, pairs $(M, D), \, |{\rm det}\, M| = |D| = m$, for which the corresponding   
 attractors are a parallelepipeds. Such pairs are all 
classified in~\cite{Z2}. For them, of course, $\bs=1$. Conjecture~\ref{conj.10}  claims that 
for all other attractors~$\bs < 1$. This means, in particular,  that 
an attractor which is not a parallelepiped cannot be presented as a finite union of 
regular sets (either convex or with a piecewise-smooth boundary). 

We can prove only the univariate version of Conjecture~\ref{conj.10}:

\begin{theorem}\label{th.50}
If an  attractor $G \subset \re$ is such that~$\bs = 1$, then~$G$ is a segment.   
\end{theorem}
We prove more: if $\bal(G) = 1$, then~$G$ is a segment. The proof uses some  facts from the theory of univariate refinement equations, from approximation theory, and from combinatorics.  Since $M$ is a number, we assume in the proof that~$M > 0$, 
the case of negative~$M$ is considered in the same way. Thus, $M = m$. 
We begin with proving several auxiliary results.

\begin{lemma}\label{l.30}
Let a set $P \subset \re$ consist of finitely many disjoint segments
with integer ends. 
Suppose several translates of~$P$ form a disjoint (up to sets of measure zero)
partition of  some segment. 
Then all the segments of the set~$P$ have the same length and 
all distances between them are multiples of that length.   
\end{lemma}
In the proof it will be convenient to use 
words ``left'' and ``right'' for the standard orientation on the real line. 

{\tt Proof}. Without loss of generality it can be assumed that 
translates of~$P$ to positive numbers cover a segment without overlaps. 
   Denote the most left segment of~$P$ by $\alpha$
and the next segment by~$\beta$. 
The distance between $\alpha$ and $\beta$ must be filled with several translates of 
$\alpha$, hence this distance is a multiple of $|\alpha|$. 
The first translate of~$P$ maps the segment $\beta$ to the segment 
$\beta +|\alpha|$. The gap between those two segments is of length $|\alpha| - |\beta|$, 
it must be filled  with several translates of the segment~$\alpha$. 
Therefore, $|\alpha| - |\beta| = k|\alpha|$ for some non-negative integer~$k$. 
Hence~$k=0$ and $|\beta| = |\alpha|$. Then by the same argument we show that the next 
segment of $P$ has length~$|\alpha|$ and that the distance from~$\beta$ to that segment is a multiple of~$|\alpha|$, etc. 
 
{\hfill $\Box$}
\smallskip

A compactly-supported function $f: \re \to \re$ is said to {\em satisfy the Strang-Fix condition of order~$\ell \in \z\cup \{0\}$} if linear combinations of its integer 
translates generate all algebraic polynomials of order~$\le \ell$. 
This condition is important in approximation theory, see~\cite{SF}.   
\begin{lemma}\label{l.40}
 A characteristic function of a compact set cannot satisfy the Strang-Fix condition
 of order bigger than zero. 
\end{lemma}
{\tt Proof}. We need to show that integer translates of~$f$ cannot generate 
a linear function. Assume the contrary: 
$\sum_{k\in \z}a_kf(x-k)\, \equiv \, x$ for some coefficients~$\{a_k\}_{k\in \z}$. Let  ${\rm supp} \,f \, \subset \, [-N,N]$. 
Then for every  $x\in [0,1]$ we have: $\sum_{k=-N}^N a_kf(x-k)\, = \, x$. 
On the other hand,~$f$ takes only the values $0$ and $1$, hence the function
$\sum_{k=-N}^N a_kf(x-k)$ takes finitely many values. This is a contradiction since 
the function $x$ takes infinitely many values on~$[0,1]$.

{\hfill $\Box$}
\smallskip 

The following result is crucial in the proof of Theorem~\ref{th.50}. 

\begin{prop}\label{p.40}
If $\bal(G) = 1$, then $G$ is a union of segments with integer ends. 
\end{prop}
{\tt Proof}. By \cite[Theorem~7.1.1]{NPS}, if a refinable function~$\varphi$
is such that $\bal_1 (\varphi) = 1$ and it does not satisfy 
the Strang-Fix condition of order~$1$, then it is $L_1$-Lipschitz, i.e  there is a 
constant $C > 0$ for which 
\begin{equation}\label{eq.Lip1}
\|\varphi(x+h) \, - \, \varphi(h)\|_1\ \le \ C\, h\, 
, \quad h \ge 0.
\end{equation}
 By Lemma~\ref{l.40}, this is true for the function $\varphi = \chi_{G}$.
 We take the system of digits~$\Delta = \{0,\ldots , M-1\}$
 with the corresponding tile~$Q = [0,1]$. Then, up to an integer translate, it can be assumed that all digits from~$D$ are nonnegative and as always~$d_0 = 0$. 
 Let $N$ be the number bigger by one than the largest digit form~$D$. 
 Since, $D \subset [0, N-1]$, we have $G \, = \, {\rm supp} \, \varphi \, \subset \, 
 [0,N-1]$. Since $\cup_{k = 0, \ldots , N-1}(Q+k)\, = \, [0,N]$, 
 we see that the set~$S = \{0, \ldots , N-1\}$ is admissible. 
We consider the corresponding vector-function~$v(x) \, = \, \bigl(\varphi(x), \ldots ,  \varphi(x+N-1)\bigr)^T \, \in \, \re^N$, and the transition 
$N\times N$ matrices $T_0, \ldots , T_{M-1}$
defined in formula~(\ref{eq.T}) for the sequence~$c_k = 1$ of $k \in D$
and $c_k = 0$ otherwise. For every $h \in (0,1)$, denote 
$\varphi_h (x) \, = \, \varphi (x+h)\, - \, \varphi (x)$
and $v_h(x) \, = \, \bigl(\varphi_h(x), \ldots ,  \varphi_h(x+N-1)\bigr)^T $. 
Applying the self-similarity equation~(\ref{eq.ss1}) we obtain 
\begin{equation}\label{eq.ssh1}
v_{M^{-1}h}(x)  \ = \ T_{\delta} v_h (Mx - \delta)\,  , \quad \delta \in \Delta\, .
\end{equation}
Denote by $\cR_k$ the set of all  $M$-adic rational numbers of 
order $k$ on the interval $[0,1)$, thus 
$R_k \, = \, \{0.\delta_1\ldots \delta_{k}\ | \ \delta_i \in \Delta, \, 
i = 1, \ldots , k\}$. For every $q \, = 0.\delta_1\ldots \delta_{k} \, \in \, \cR_k$, let $\Pi_q \, = \, 
T_{\delta_1}\cdots T_{\delta_k}$. 
Iterating $k$ times equation~(\ref{eq.ssh1}) we obtain  
\begin{equation}\label{eq.ssh2}
v_{M^{-k}h}(x) \ = \ \Pi_q \, v_h\bigl(M^k(x - q)\bigr) \, , 
\qquad x \in \bigl[q, q+M^{k} \bigr)\, , \ q \in \cR_k\, . 
\end{equation}
Denote $A_{\delta} = T_{\delta}|_{U}$. It was proved 
in~\cite[Theorem 5.2.2]{NPS}  that 
a refinable function~$\varphi$ is ~$L_1$-Lipschitz if and only if 
for every $k$, we have $\sum_{q \in \cR_k}\|\Pi_q|_{U}\|\, \le \, C_0$, where $C_0$ is some constant. 
On the other hand, all the products $\Pi_q, \, q\in \cR_k$, are simple matrices. 
Denote by $C_1$ the smallest positive norm of all simple matrices restricted 
to~$U$. Thus, 
 for every $q \in \cR_k$, either $\Pi_q|_{U} = 0$ or $\|\Pi_q|_{U}\| \ge C_1 $. 
 Hence, for every $k\in \n$, the set $\{\Pi_q|_{U}\ | \ q \in \cR_k\}$ contains at most 
 $r \, = \, \bigl[\frac{C_0}{C_1}\bigr]$ nonzero operators. 
Apply this fact to equation~(\ref{eq.ssh2}) taking into account that 
the vector $v_h\bigl(M^k(x - q)\bigr)$ belongs to $U$. We see that on all 
but $r$ intervals $[q, q+M^{-k})$, the function $v_{M^{-k}h}(x)$
is an identical zero. This holds for every~$h \in (0,1)$ and the set of segments 
on which this function is zero is the same for all~$h$. Therefore 
on all but $r$ those segments, the function $v(x)$ is an identical constant. 
This is true for all~$k\in \n$. Hence, the function~$v(x)$ is piecewise-constant 
with at most $r$ points of discontinuity.   
 Consequently, the function~$\varphi$ is piecewise-constant with at 
 most $r+N+1$ points of discontinuity (in the $N+1$ integer points 
 of the segment~$[0,N]$ the function $\varphi$ may also be discontinuous).    
 By \cite[Theorem~1]{LLS} for every piecewise-constant refinable function with
 finitely many points of discontinuity, all those points are integer.   
 Hence, $G$ is a union of segments with integer ends.

{\hfill $\Box$}
\smallskip

{\tt Proof of Theorem~\ref{th.50}}. By Proposition~\ref{p.40},
 the set~$G$ consists of several segments 
with integer ends.  Without loss of generality it can be assumed that 
the number $M$ is positive and exceeds the diameter of the 
set~$G$. Otherwise we 
iterate the refinement equation for $\varphi = \chi_{G}$ several times, say $k$ times, and 
obtain a refinement equation with the factor $M^k$, then we just replace $M$
by $M^k$. Denote by $\alpha$ the most left segment of~$G$. 
We have $G = \cup_{d\in D}M^{-1}(G+d)$, hence 
$MG \, = \, \cup_{d\in D}(G+d)$. The most left segment of 
$G$ is $M\alpha$, its  length~$M|\alpha|$ exceeds the diameter of~$G$
(since $M > {\rm diam}\, (G)$). 
Moreover, the distances from this segment to other segments 
of~$MG$ being multiples of the number~$M$ also 
exceed the diameter of~$G$. Therefore, $M\alpha \, = \, \cup_{d\in D'}(G+d)$, 
where $D'$ is some subset of~$D$. Thus, several translates of the set~$G$
form the segment~$M\alpha$. Applying now Lemma~\ref{l.30} to the set~$P = G$, we 
obtain that  all  segments of the set~$G$ have the same length~$|\alpha|$ and all distances between them are multiples of this number.
Therefore, all distances between the segments of the 
set~$MG$ are multiples of~$M|\alpha|$. Hence, each segment of the set 
$MG$ is $M\alpha + kM|\alpha|$ with some $k\in \n$. Consequently, 
for that segment we have $M\alpha + kM|\alpha| \, = \, 
\cup_{d\in D'+kM|\alpha|}(G+d)$. Therefore, $D'+kM|\alpha| \, \subset \, D$. 
However, the sets $D'$ and $D'+kM|\alpha|$ are equal modulo 
$M$. This is impossible, since all elements of the digit set~$D$ are different 
modulo~$M$.

{\hfill $\Box$}
\smallskip

\begin{center}
\textbf{8. Application to synchronising automata}
\end{center}
\bigskip 

The theory of synchronising automata originated in 1960s   has found 
numerous applications in engineering and computer science. 
It is actively developing in the modern literature, see~\cite{CHJ, K, T, V}
and references therein.

Suppose some system can be at $N$ different states. 
There are $m$ actions that change the states of the system. 
The $k$th action changes the states according to a prescribed mapping~$f_k$ defined on the set of states. If we enumerate the states by numbers 
from $1$ to $N$, then the $k$th action changes the $j$th state to the state 
$f_k(j)$, $j=1, \ldots , N$. The set of states and of actions is called 
{\em deterministic finite automaton}. It can be defined by 
the directed graph $G(V,E)$ with coloured edges. The vertices correspond to the 
states and edges of each color are associated to actions. 
There is an edge of $k$th colour from the vertex $j$ to the vertex~$i$
if $f_k(j) = i$. Thus, $G$ has $N$
vertices and $m$ outgoing edges from every vertex -- one edge of each colour. 
The edges of the $k$th colour generate an adjacency matrix~$B_k$.  
We have $(B_k)_{ij} = 1$ if $f_k(j) = i$ and  $(B_k)_{ij} = 0$ if 
$f_k(j) \ne i$. 
Thus, the matrix~$B_k$ is simple in the sense of Definition~\ref{d.simple}: 
every column possesses exactly one element equal to one and all others 
being zeros. We see that a deterministic finite automaton is completely defined 
by the family of simple $N\times N$ matrices $\cB = \{B_0, \ldots , B_{m-1}\}$. 
Conversely, every family of simple  $N\times N$ matrices defines an automaton. 
Here it will be more convenient to enumarate the actions (colours) not from $1$ to $m$
as in the most of literature on automata but from $0$ to $m-1$. 

A finite sequence of actions (colours) is called a {\em synchronising sequence}
or {\em reset word} if application of this sequence of actions sends the system to 
one and the same state, independently of the initial state. In terms of the matrices~$\{B_k\}_{k=0}^{m-1}$ a reset word is a sequence of numbers $k_1, \ldots , k_s$
from $\{0, \ldots , m-1\}$
such that the corresponding product $\Pi = B_{k_s}\cdots B_{k_1}$ has a row of ones. 
Since $\Pi$ is a simple matrix it follows that all other entries
of~$\Pi$ are zeros. 

In practice a reset word  allows the user to make a reset the system i.e., sending  it to the initial state even if  its current state is unknown. There are lots of applications of this notion in computer science, electronics, robotics, etc. 
There are efficient polynomial time algorithms to decide the existence of a reset word 
and to find it~\cite{V}. On the other hand, finding the  shortest possible reset word is an NP-complete problem~\cite{E, V}. There is a famous \v{C}ern\'y conjecture (1964)
claiming that if a reset word exists then the shortest reset word has length at most 
$(N-1)^2$. This lower bound is sharp~\cite{Cer}. The conjecture is still open and the best known upper bounds is cubic in~$N$~\cite{S}.

Now come back for a moment to the self affine attractors. Let~$G$ be an attractor. The transition matrices $T_{\delta}\, , \ \delta \in \Delta$, are all simple. Therefore the family~$\cT$
of these matrices generates a deterministic finite automaton. What is the sense 
of the the surface regularity~$\bs(G)$ is terms of the automaton of the family~$\cT$? Can it be of interest to the automata theory? To answer this question we 
introduce a concept of parameter of synchronisation. 
In the following theorem we use the same subspace 
$W \, = \, \{x \in \re^N \ | \ \sum_{i=1}^N x_i = 0\}$
and denote $A_k = B_k|_W$, where $B_k$ is an adjacency matrix of an automaton. 
\begin{theorem}\label{th.60}
Let a deterministic finite automaton be given. For a natural 
$k$, let $P_k$ be the probability that a random word of length~$k$
of the alphabet~$\{0, \ldots , m-1\}$ is not a reset word. 
Then there exists a limit~$p \, = \, \lim_{k\to \infty}[P_k]^{1/k}$. 
This limit is equal to the spectral radius~$\rho(\mathbf{A})$ of the operator
\begin{equation}\label{eq.AN}
\mathbf{A}(X) \ = \ \frac{1}{m}\, \sum_{i=0}^{m-1} A^*_iXA_i\, , \qquad X \in \cM_{N-1}\, . 
\end{equation}
 acting on the $N(N-1)/2$-dimensional space~$\cM_{N-1}$ of symmetric 
 matrices of size~$N$.  
 
 The automaton has a reset word if and only if $\rho(\mathbf{A}) < 1$. 
\end{theorem}
{\tt Proof.} Denote by $C_1$ and $C_2$ respectively the smallest and the largest strictly positive norms of all simple matrices 
restricted to the subspace~$W$. Since there are finitely many simple matrices, it follows that $C>0$. If a word $\ell_1\ldots \ell_k$ is reset, then 
the product $\Pi = B_{\ell_1}\cdots B_{\ell_k}$ has a row of ones and therefore 
$\Pi|_W = A_{\ell_1}\cdots A_{\ell_k} = 0$. Otherwise, $C_1 \le \|A_{\ell_1}\cdots 
A_{\ell_k}\| \, \le C_2$. For every~$k$, the number of nonzero products among all 
products $A_{b_1}\cdots A_{b_k}$ is equal to $m^kP_k$. 
Therefore, the value 
$$
S_k \ = \ m^{-k}\sum_{b_1, \ldots, b_k} \|A_{b_1}\cdots A_{b_k}\|^2
$$
is between $C_1^2P_k$ and $C_2^2P_k$. The power $1/k$ of this value tends to
$\rho_2^2$, where $\rho_2$ is the $L_2$-spectral radius 
of the family~$A_0, \ldots, A_{m-1}$. Hence the limit $\lim_{k\to \infty}[P_k]^{1/k}$
exists and is equal to~$\rho_2^2$, which is in turn equal to the spectral 
radius of the operator~(\ref{eq.AN}). 

Since all norms $\|A_{b_1}\cdots A_{b_k}\|$ are bounded above by $C_2$, 
it follows that $\rho_2 \le 1$, and hence $p = \rho_2^2 \le 1$.
It is well-known that for every family of operators~$\{A_0, \ldots , A_{m-1}\}$, there exists a constant~$C > 0$ such that $S_k \, \ge \, C\, \rho_2^{2k}$ for every $k\in \n$. Hence if $p=1$, then $S_k \ge C$ for all $k$. However, if there exists at least one zero product of those operators, then $S_k \to 0$ as $k \to \infty$. 
Therefore, if $p=1$, then there is no zero product, which means that 
there is no product of operators $B_0, \ldots , B_{m-1}$ with a row of ones, 
i.e., there is no reset word.

{\hfill $\Box$}
\smallskip

We call the number $p$ the {\em parameter of synchronisation}
of the automaton. It has the following meaning. Assume we do not know a reset word and instead take a  random sequence of actions of length~$k$;  
 then  we obtain a reset word apart from the probability approximately~$\, p^k$.  Thus, the parameter $p$ shows the 
degree of random synchronisation of an automaton.  As we see from Theorem~\ref{th.60}, this parameter can be effectively computed merely by finding the largest 
eigenvalue of operator~(\ref{eq.AN}). The following theorem reveals a curious relation between the parameter of synchronisation and the surface regularity of 
a self-affine tile.

\begin{theorem}\label{th.70}
Let a  tile $G(M, D)$ be given and its dilation matrix
$M$ be isotropic; then  
for the automaton defined by the transition matrices~$T_{\delta}, \, 
\delta \in \Delta$, the parameter of synchronisation~$p$
is equal to $r^{- \bs}$, where $r = \rho(M)$ and $\bs$ is the surface regularity 
of~$G$.  
\end{theorem}
{\tt Proof}. Applying Theorem~\ref{th.tile-isotr} 
we obtain $\bs = 2\log_{1/r} \rho_2$ and hence $\rho_2^2 \, = \, r^{- \bs}$. 
By Theorem~\ref{th.60}, $p = \rho_2^2$, which completes the proof.  

{\hfill $\Box$}
\smallskip 

 \textbf{Acknowledgements}. The author is grateful to T.Zaitseva for her help in computations and in illustrations and to R.Karasev and S.Ivanov for useful discussions. 
 We are thankful to the anonymous Referee for his attentive reading and valuable comments. 
 
  A part of the work was done during the visit of the 
author in the Erwin Schr\"odinger Institute (ESI), Vienna, Austria. 
We express our thanks to the institute for hospitality.

\end{document}